\documentclass[12pt]{article}
\usepackage[T1]{fontenc}
\usepackage[left=25mm,right=25mm,top=25mm, bottom=25mm]{geometry}
\usepackage{amssymb,amsmath,amscd, amsthm,stmaryrd,verbatim,
	mathrsfs,tikz,color , amsfonts}
\usepackage{cite}
\usepackage{setspace} 
\allowdisplaybreaks
\usepackage[vcentermath,enableskew,noautoscale]{youngtab}
\usepackage{young}
\usepackage{hyperref}
\usepackage{pifont}

\newtheorem{lem}{Lemma}[section]
\newtheorem{them}[lem]{Theorem}

\newtheorem{cor}[lem]{Corollary}

\title{\textbf{Periodic points in the $\beta$-transformation with a hole at 0}}
\author{ Yuzheng Bi 
	\\
\small\textit{College of Mathematics and Statistics, Chongqing University, Chongqing 401331, China}}
\date{}
\begin{document}
	\bibliographystyle{plain}
	\maketitle
	\renewcommand{\abstractname}{}
	\begin{abstract}
		\noindent\textbf{Abstract:}	
		For $\beta\in(1,2]$ let $T_\beta: [0,1)\to[0,1); x\mapsto \beta x\pmod 1$. In this paper we study the periodic points in the open dynamical system $([0,1), T_\beta)$ with a hole $[0,t)$. For $p\in\mathbb{N}$ we characterize the largest $t$, denoted by $S_\beta(p)$, in which the survivor set $K_\beta(t)$ has a periodic point of smallest period $p$. More precisely, we give precise formulae for this critical value $S_\beta(p)$ when $\beta=2$, $\beta=\frac{1+\sqrt{5}}{2}$ and $\beta$ being the tribonacci number. We show that for $\beta=2$ the critical value $S_2(p)$ converges to $1/2$ as $p\to \infty$.  When $\beta=\frac{1+\sqrt{5}}{2}$, the critical value $S_\beta(p)\to \frac{1}{\beta^3-\beta}$. While $\beta$ is the tribinacci number, the critical value $S_\beta(p)\to \frac{\beta^{2}+1}{\beta^4-\beta}$.

	\end{abstract}
	
\section{Introduction}
{The study of numbers has a long history. Initially, people used binary, octal, decimal, and hexadecimal systems for research. These different expressions share a common feature that the basis used for expanding real numbers is an integer greater than 1.\par
The expansion of real numbers on non integer bases was first studied by Rényi\cite{AR57}, for non-integer $\beta$, any real number $x$ can be represented as
\begin{equation} \label{introduction}
	x=\epsilon_0(x)+\frac{\epsilon_1(x)}{\beta}+\frac{\epsilon_2(x)}{\beta^2}+\cdots+\frac{\epsilon_n(x)}{\beta^n}+\cdots
\end{equation}
The numbers can be obtained by the following recursive formula
$$\epsilon_0(x)=[x],\quad r_0(x)=(x)$$
$$\epsilon_i(x)=[\beta\cdot r_{i-1}(x)],\quad r_i(x)=(\beta\cdot r_{i-1}(x)),\quad\left\{i=1,2\cdots\right\}$$ 	where $[x]$ represents the integer part of $x $, $(x)$ represents the decimal part of $x$. (For more information of the properties of $\beta$ expansion, please refer to the literature\cite{BY99,HF78}.) 

Specially, $\forall\beta>1$, let $T_\beta:[0,1)\rightarrow[0,1)$  is the $\beta$-transformation\cite{BR15}, defined by:
\begin{equation}
	T_\beta(x)=\beta x-[\beta x],
\end{equation}
We obtained a dynamical system, and by iterating this mapping, we can obtain the $\beta$ expansion of $x$:
\begin{equation}
	x=\sum\limits_{i=1}^\infty\frac{x_i}{\beta^i},
\end{equation}
where $x_i=[\beta T_\beta^{i-1}(x)]$, and$x_i\in \left\{0,1,\cdots[\beta-1]\right\}$.\par
Afterwards, Rényi proved that $T_\beta$ is ergodic with respect to the Lebesgue measure m, so that $\forall f\in L^1[0,1)$, and $m-a.e.x$, we have 
\begin{equation}
	\lim\limits_{n\rightarrow \infty}\frac{1}{n}\sum\limits_{i=0}^{n-1}f(T_\beta^{i-1}(x))=\int f\ d\nu,
\end{equation}
where $\nu$ is the Parry measure equivalent to the Lebesgue measure. Parry\cite{PW60} further characterized the traversal measure in $\nu$ 1960. (For more traversal theorems, see reference \cite{BK31, BI65}.)\par

Since this non-integer basis was introduced, this real number expansion based on non-integer basis has been widely studied from different perspectives, which are related to fields such as symbolic dynamical systems, fractal geometry, number theory and so on.\cite{IY74,KV07,SD15,FT12,CB92,KS03}.\par
The research on the unique $\beta$-expansion has always been a focus research, many scholars have studied the topological structure, metric properties, and fractal dimension of real number sets with unique expansion.\cite{MV1112,KV1113,KP9814,DW1515,MV0916}.\par
Allouche, Jean-paul and Clarke\cite{JM09} studied the periodic point problem with unique $\beta\in(1,2]$ expansion, by studying any given positive integer $n\geq2$, they obtain the smallest $\beta_n\in [\frac{1+\sqrt{5}}{2},2)$, such that $\forall \beta>\beta_n$ there exists at least one unique $\beta$ expansion with the smallest period $n$ in 2009. Tan and Ge\cite{TG17} promoted the work of reference[19] in 2017, for any $\beta>1$, they studied the periodic problem with the unique $\beta$ expansion. Urba\'{n}ski first studied the $\beta$ expansion problem of digging a hole at the origin in his work \cite{UR86, UR87} in 1986 and 1987. Kalle, Kong, Langeveld and Li\cite{CK20} studied the $\beta$ expansion with digging a hole $[0,t)$ on $[0,1)$ in 2020. They studied the set-valued bifurcation set. Baker and Kong\cite{KD20} studied the dimension bifurcation set which is related  to dynamical systems in 2020.\par 
Based on previous research, this article focuses on studying the value of $S_\beta(p)$ which is defined as following:\par 
For $\beta\in (1,2]$ the $\beta-$transformation $T_\beta : [0,1) \rightarrow [0,1)$ is defined by $T_\beta(x)=\beta \cdot x(mod 1)$. \\
For $t\in[0,1)$, let $K_\beta(t)$ be the survivor set of $T_\beta$ with a hole $(0,t)$ which is given by 
$$K_\beta (t):=\left\{x \in [0,1):T^n_\beta(x) \geq t,\forall n \geq 0\right\}.$$
For any given $p\in \mathbb{N}^+$, let $S_\beta(p)$ be given by
$$S_\beta(p)=\sup\left\{t\in[0,1):K_\beta(t) \mbox{contains a point of the smallest period $p$}\right\}.$$
We obtained the value of $S_\beta(p)$ for any positive integer $p$ when $\beta=2$ in Theorem 3.1. When $\beta=\frac{1+\sqrt{5}}{2}$, we finally obtained the value of $S_\beta(p)$ for any positive integer $p$ in Theorem 4.3. When $\beta$ is the tribonacci number, we finally obtained the value of $S_\beta(p)$ for any positive integer $p$ in Theorem 5.8. \par 
In the next section, we will introduce some basic definitions and preliminaries which is related to the research of this paper. In the third section, we show that for $\beta=2$ the critical value $S_2(p)$ converges to $1/2$ as $p\to \infty$. In the fourth section, we show the critical value $S_\beta(p)\to \frac{1}{\beta^3-\beta}\approx 0.38212$ as $p\to \infty$ when $\beta=\frac{1+\sqrt{5}}{2}$. In the fifth section, we show the critical value $S_\beta(p)\to \frac{\beta^{2}+1}{\beta^4-\beta}\approx 0.45626$ as $p\to \infty$ when $\beta$ is the tribonacci number. 

}

 \section{Basic definitions and Preliminaries}
Let $\left\{0,1\right\}^\mathbb{N}$ be the set of sequences composed of digits 0 and 1, and define the left-shift operator on $\left\{0,1\right\}^\mathbb{N}$$\sigma$ by $$\sigma((x_i))=(x_{i+1}).$$\par
In this paper, for any two sequences $(x_i),(y_i)\in\left\{0,1\right\}^\mathbb{N}$, we define the lexicographic order symbol $\prec$ between sequences: if $x_1< y_1$ or there exists $m \in \mathbb{N}_\geq2$, satisfying $\forall j \leq m,x_j=y_j$ and $x_{j+1}<y_{j+1}$, we say $(x_i)\prec(y_i)$. Similarly, we define $"\succ"$.\par
In addition, we refer to $(x_i)\preceq(y_i)$, if$(x_i)\prec(y_i)$ or $(x_i)=(y_i)$. Similarly, we define $(x_i)\succeq(y_i)$, if$(x_i)\succ(y_i)$ or $(x_i)=(y_i)$.\par
We denote the greedy $\beta$ expansion of $x$\cite{KV12}by $b(x,\beta)$, which is obtained by always selecting the maximum possible value at each step $b_i(x,\beta)$. The greedy $\beta$ expansion of $x$ is the largest $\beta$ expansion in the sense of lexicographic order.\par 
We denote the quasi-greedy $\beta$ expansion of $x$ by\cite{DK95}$\widetilde{b}(x,\beta)$, which is the largest infinite $\beta$ expansion in the sense of lexicographic order which does not end with $0^\infty$.\par
The quasi-greedy $\beta$ expansion of 1 plays the important role in many researches. We denote the quasi-greedy $\beta$ expansion of 1 by $\delta(\beta)=\delta_1(\beta)\delta_2(\beta)\cdots$, and we point out that when the greedy $\beta$ expansion of $x$ is $b(x,\beta)=b_1b_2\cdots b_n0^\infty$, the corresponding quasi-greedy $\beta$ expansion of $x$ is $\widetilde{b}(x,\beta)=b_1b_2\cdots b_n^-\delta(\beta)$, where $b_n^{-}=b_n-1$.\par
For any $\beta>1$, suppose that the $\beta$ expansion of real number $x$ is
\begin{equation}
	x=\sum\limits_{i=1}^\infty\frac{x_i}{\beta^i},
\end{equation}
We record it in this paper as:
\begin{equation}
	x=\sum\limits_{i=1}^\infty\frac{x_i}{\beta^i}=:[(x_1x_2x_3\cdots)]_\beta.
\end{equation}

Multinacci number is an extended form of Fibonacci sequence, which belongs to a kind of multi order recursive sequence, which is the unique root in $(1,2)$ of the equation
$$x^{m+1}=x^m+x^{m-1}+\cdots+x+1.$$
multinacci number has any order form $m(m\in\mathbb{N}_{>1})$: the second order multinacci number is the gold ratio, the third order multinacci number is Tribonacci number\cite{CJ15,WR82} and  the fourth multinacci number is Tetranacci number\cite{GD2329}.\par
Importantly, if $\beta$ is a multinacci
number then its quasi-greedy expansion of 1 is of the form $((1^m0))^\infty$. This property will be
useful in our analysis. For simplicity of description, we denote the golden ratio as $\phi_2$ and denote the tribinacci number as $\phi_3$ in this paper.\par

\section{The case for $\beta=2$}
The real number expansion based on 2 is a common problem in the study of $\beta$-expansion. The transformation involved here is the $T_2$ transformation. In this chapter, we will study the value of $S_2(p)$ which is related to the survivor set $K_2(t)$ at digging a hole $[0,t)$ at the origin. For more research related to the survival set after digging a hole at origin under the $T_2$ transformation, please refer to the literature\cite{BY1130,DET1331,AB1432,SID1433,GS1534,CT1735}.

\begin{them} \label{th:01}
	For $\beta=2,S_2(p)=[(01^{p-1})^\infty]_2$, for any positive integer $p$.
	\begin{proof}
		{
			\begin{itemize}
				\item[(1)]
				{Firstly, we will prove that $S_2(p)\geq [(01^{p-1})^\infty]_2=: t_*$.\\
					We get that  $$\widetilde{K_\beta(t_*)}=\left\{(d_i)\in \left\{0,1\right\}^\mathbb{N} :(01^{p-1})^\infty\leq\sigma^n((d_i))\right\}.$$
					Obviously the sequence $(01^{p-1})^\infty\in\widetilde{K_\beta(t_*)}$. That is $\widetilde{K_\beta(t_*)}$ contains a sequence of the smallest period p.\\
					By the definition of $S_\beta(p)$, we have $S_\beta(p)\geq t_*.$
				}
				
				\item[(2)] Next to prove that $S_2(p)\leq [(01^{p-1})^\infty]_2=: t_*$. \par 
				Suppose on the contrary that $S_2(p)>t_*$.
				Then there exists $t\in (t_*,S_2(p))$, and then $\widetilde{K_2(t)}$ contains at least a periodic sequence $(\omega_i)$ of the smallest period p. And we have 
				\begin{equation}
					(01^{p-1})^\infty\prec(t_i)\preceq\sigma^j((\omega_i))\prec(1)^\infty, \forall j \in \mathbb{N}.
				\end{equation}
				\begin{itemize} 
					\item [Case 1:]{
						If $(\omega_i)=(\omega_1 \omega_2 \cdots \omega_p)^\infty$ with $\omega_1=1$ and $\omega_i\in \left\{0,1\right\}$ for $i\neq 1$, there exists $k\in\left\{1,2\cdots p\right\}$ such that $\omega_k=0$. Then
						$$T_2^{k-1}(\omega_i)=(0\widehat{\omega_2}\widehat{\omega_3}\cdots\widehat{\omega_p}).$$
						Hence we have $T_2^{k-1}(\omega_i)=(0\ \widehat{\omega_2}\ \widehat{\omega_3}\cdots\widehat{\omega_p})\preceq (01^{p-1})^\infty$, which is contradicting to the inequality (7).
					}
					\item [Case 2:]{
						If $(\omega_i)=(\omega_1 \omega_2 \cdots \omega_p)^\infty$ with $\omega_1=0$, and $\omega_i \in \left\{0,1\right\}$ for $i\neq 1$, we have that $(\omega_i)= (0\omega_2 \omega_3\cdots \omega_p)^\infty\preceq(01^{p-1})^\infty$ which is contradicting to the inequality (2.1).  
					}
				\end{itemize}
			\end{itemize}
			Thus we have $S_2(p)\leq t_*$.\par
			To sum up, we get $S_2(p)=[(01^{p-1})^\infty]_2$, for any positive integer p.
		}
	\end{proof}  
\end{them}

{By calculation, we get that $$S_2(p)=\frac{2^{p-1}-1}{2^p-1},\forall p\in\mathbb{N}_{\geq1},$$ which is increasing to $\frac{1}{2}$ as p increases to infinity.}

\begin{figure}[htb]  
	\centering  
	\begin{minipage}{.6\textwidth}  
		\centering  
		\includegraphics[width=.8\linewidth]{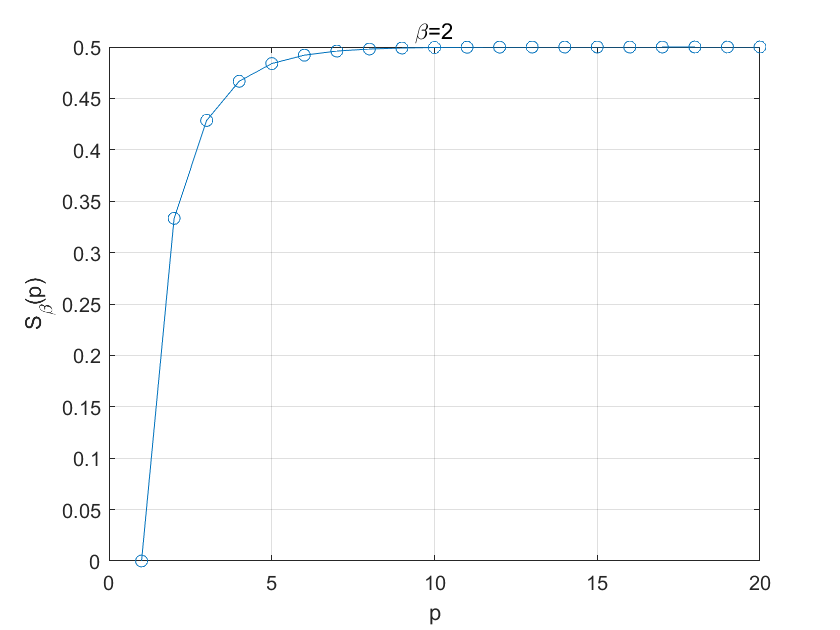}  
		\label{fig:image2}  
	\end{minipage}
	\tiny\caption{The value of $S_\beta(p)$ for $\beta=2$ with the period p from 1 to 20.}
\end{figure}

\section{The case for $\beta=\phi_2$ }\label{sec 4}
Next we will consider the value of $S_\beta(p)$ when $\beta$ is the golden number. Here the golden ratio $\beta \approx 1.61803$is the real root of the polynomial $x^2-x-1=0$. Since the quasi-greedy $\beta$-expansion of 1 is $(10)^\infty$ for $\beta=\frac{1+\sqrt{5}}{2}$, we get
$$\widetilde{K_\beta(t)}=\left\{(d_i)\in \left\{0,1\right\}^\mathbb{N} :(t_i)\preceq\sigma^n((d_i))\prec(10)^\infty\right\},$$
where $\sigma$ is the left-shift operator.

\begin{lem}
	For $\beta=\phi_2$, we have $$S_\beta(p)=[(0(01)^m)^\infty]_\beta ,p=2m+1,m\in \mathbb{N}.$$
	\begin{proof}
		\begin{itemize}
			\item [(1)]Since there are only two 1-periodic sequences $(0)^\infty$ and $(1)^\infty$ which are only composed by 0 and 1 and due to $(1)^\infty\succ(10)^\infty$, we get $S_\beta(1)=0$ immediately. 
			\item [(2)]Then we will prove that $S_\beta(2m+1)\geq [(0(01)^m)^\infty]_\beta=: t_*,p=2m+1,m\in \mathbb{N}_{\geq1}$.\\
			We get that  $$\widetilde{K_\beta(t_*)}=\left\{(d_i)\in \left\{0,1\right\}^\mathbb{N} :(0(01)^m)^\infty\preceq\sigma^n((d_i))\prec(10)^\infty\right\}.$$
			Obviously the sequence $(0(01)^m)^\infty\in\widetilde{K_\beta(t_*)}$. That is $\widetilde{K_\beta(t_*)}$ contains a sequence of the smallest period 2m+1.\\
			By the definition of $S_\beta(2m+1)$ we have $S_\beta(2m+1)\geq t_*.$ \\ 
			Finally we will prove that $S_\beta(2m+1)\leq [(0(01)^m)^\infty]_\beta=: t_*,m\in \mathbb{N}_{\geq1}$.\\
			Suppose on the contrary that $S_\beta(2m+1)>t_*$, then there exists $t\in(t_*,S_\beta(2m+1))$ and $\widetilde{K_\beta(t)}$ contains at least a periodic sequence $(\omega_i)$ of the smallest period 2m+1. We have
			\begin{equation}
				(0(01)^m)^\infty\prec(t_i)\preceq \sigma^j((\omega_i))\prec(10)^\infty,\forall j \in \mathbb{N}.
			\end{equation}
			Since $(\omega_i)$ is of the smallest period $2m+1,m\in \mathbb{N}_{\geq1}$, there is at least one 1 among $\omega_1,\omega_2\cdots\omega_{2m+1}$.
			\begin{itemize} 
				\item[Case 1:] {
					If $(\omega_i)=(\omega_1 \omega_2 \cdots \omega_{2m+1})^\infty$ with $\omega_1=1$ and $\omega_i\in \left\{0,1\right\}$ for $i\neq 1$, we can write
					$$(\omega_i)=(\omega_1 \omega_2 \cdots \omega_{2m+1})^\infty=(1\ \square\ 1\ \square\cdots1\ \square)^\infty,$$
					where each $\square$ only contains 0. Then $\exists l \in \mathbb{N}_{\geq1}$ such that $$\sigma^l((\omega_i))=(\square\ 1\ \square\ 1\cdots\ \square\ 1)^\infty.$$
					Since $\sigma^j((\omega_i))\prec(10)^\infty,\forall j \in \mathbb{N}$, there exists at least one 0 in each $\square$. And if each $\square$ contains exactly one 0, the sequence $\sigma^l((\omega_i))=(01)^\infty$ which is a contraction to $(\omega_i)$ is of the smallest period 2m+1. Therefore there exists a $\square$ contains more than one 0. Then $\exists k \in \mathbb{N}$ such that $\sigma^{l+k}((\omega_i))=(00\widetilde{\omega_3}\cdots\widetilde{\omega_{2m+1}})\preceq (0(01)^m)^\infty$ which is a contraction to the inequality (8).
				}
				\item [Case 2:]{
					If $(\omega_i)=(\omega_1 \omega_2 \cdots \omega_{2m+1})^\infty$ with $\omega_1=0$, and $\omega_i \in \left\{0,1\right\}$ for $i\neq 1$, we can write 
					$$(\omega_i)=(\square\ 1\cdots \square\ 1\ \square)^\infty.$$  
					The same to the Case 1, we get a contraction.
				}
			\end{itemize}
			Hence we have $S_\beta(2m+1)=[(0(01)^m)^\infty]_\beta,\forall m\in \mathbb{N}.$
		\end{itemize}			
	\end{proof}
\end{lem}

\begin{lem}
	For $\beta=\phi_2$ , we have $$S_\beta(p)=[(001(01)^{m-2}001(01)^m)^\infty]_\beta ,p=4m+2,m\in \mathbb{N}_{\geq2}.$$
	\begin{proof}
		\begin{itemize}
			\item[(1)]
			{$S_\beta(2)=0.$\\
				Firstly, $S_\beta(2)\geq 0$ is trivial. Next to prove that $S_2(p)\leq 0$.\par Suppose on the contrary that $S_2(p)>0$.\\
				Then there exists $t\in (0,S_2(p))$, and then $\widetilde{K_\beta(t_*)}$ contains at least a periodic sequence $(\omega_i)$ of the smallest period 2. And we have 
				$$(0)^\infty\prec(t_i)\preceq \sigma^j((\omega_i)).$$
				There are only $(10)^\infty$ and $(01)^\infty$ two periodic sequences of the smallest period 2. Hence $(\omega_i)=(10)^\infty$ or $(\omega_i)=(01)^\infty$ which are contradict to $\sigma^j((\omega_i))\prec(10)^\infty,\forall j \in \mathbb{N}$.
				Therefore we have $S_2(p)\leq 0$.
				Furthermore we have $S_2(p)= 0.$			 		
			}
			
			\item[(2)]
			{$S_\beta(6)=[(00(01)^2)^\infty]_\beta.$\\
				Firstly, we will prove that $S_\beta(6)\geq [(00(01)^2)^\infty]_\beta=: t_*$.\\
				We get that  $$\widetilde{K_\beta(t_*)}=\left\{(d_i)\in \left\{0,1\right\}^\mathbb{N} :(00(01)^2)^\infty\preceq\sigma^n((d_i))\prec(10)^\infty\right\}.$$
				Obviously the sequence $(00(01)^2)^\infty\in\widetilde{K_\beta(t_*)}$. That is $\widetilde{K_\beta(t_*)}$ contains a sequence of the smallest period 6.\\
				By the definition of $S_\beta(6)$ we have $S_\beta(6)\geq t_*.$\\
				Next to prove that $S_\beta(6)\leq [(00(01)^2)^\infty]_\beta=: t_*$. Suppose on the contrary that $S_\beta(6)>t_*$.
				Then there exists $t\in (t_*,S_\beta(6))$, and then $\widetilde{K_\beta(t)}$ contains at least a periodic sequence $(\omega_i)$ of the smallest period 6. And we have 
				\begin{equation}
					(00(01)^2)^\infty\prec(t_i)\preceq \sigma^j((\omega_i))\prec(10)^\infty.
				\end{equation}
				Since $(\omega_i)$ is of the smallest period 6, there is at least one 1 among $\omega_1,\omega_2\cdots\omega_6$.
				And since $\sigma^j((\omega_i))\prec(10)^\infty,\forall j \in \mathbb{N}$, there is no consecutive 1 in $(\omega_i)$.
				Due to a fact that if there are more than there 1 among $\omega_1,\omega_2\cdots\omega_6$ there will exist consecutive 1 in $(\omega_i)$, we have 
				$$\max\# \left\{i:\omega_i=1,i\in\left\{1,2,\cdots6\right\}\right\}=3.$$
				If there are three 1 among $\omega_1,\omega_2\cdots\omega_6$, due to $\sigma^j((\omega_i))\prec(10)^\infty,\forall j \in \mathbb{N},$\\
				we have $(\omega_i)=(101010)^\infty$ or $(\omega_i)=(010101)^\infty$ which are of the smallest period 2 which is a contraction.
				Hence there are less than three 1 among $\omega_1,\omega_2\cdots\omega_6$.
				If there is one 1 among $\omega_1,\omega_2\cdots\omega_6$, there exist $k\in\mathbb{N}$ such that $$\sigma^k((\omega_i))=(000001)^\infty\prec(00(01)^2)^\infty,$$ which is a contraction.\\
				Hence there is exactly two 1 between $\omega_1,\omega_2\cdots\omega_6$.We can write
				$$(\omega_i)=(\square\ 1\ \square\ 1\ \square)^\infty,$$
				where each $\square$ only contains 0.
				There exists $k\in\mathbb{N}$ such that $$\sigma^k((\omega_i))=(\square_1\ 1\ \square_2\ 1)^\infty.$$
				since $\sigma^j((\omega_i))\prec(10)^\infty,\forall j \in \mathbb{N}$, each $\square$ contains 0.
				\begin{itemize}
					\item[$(\textrm{i})$]
					{
						If there is one 0 in $\square_1$ and three 0 in $\square_2$, $\sigma^{k+2}(\omega_i)=(000101)^\infty$ which is a contraction to (9).
					}
					\item[$(\textrm{ii})$]
					{
						If there is two 0 in $\square_1$ and two 0 in $\square_2$, $\sigma^{k}(\omega_i)=(001001)^\infty$which is a sequence of the smallest period 3 which is a contraction to (9).
					}
					\item[$(\textrm{iii})$]
					{
						If there is three 0 in $\square_1$ and one 0 in $\square_2$, $\sigma^{k}(\omega_i)=(010001)^\infty$ which is a contraction to (9).
					}
				\end{itemize}
				Hence there is no such 6-period sequence $(\omega_i)$ contained in  $\widetilde{K_\beta(t)}$ which is a contraction.
				Therefore we get $S_\beta(6)=[(00(01)^2)^\infty]_\beta.$
			}
			\item [(3)]For $p=4m+2,m\in \mathbb{N}_{\geq2},S_\beta(4m+2)=[(001(01)^{m-2}001(01)^m)^\infty]_\beta.$\\
			1)Firstly, we will prove that $S_\beta(4m+2)\geq [(001(01)^{m-2}001(01)^m)^\infty]_\beta=: t_*$.\\
			We get that  $$\widetilde{K_\beta(t_*)}=\left\{(d_i)\in \left\{0,1\right\}^\mathbb{N} :(001(01)^{m-2}001(01)^m)^\infty\preceq\sigma^n((d_i))\prec(10)^\infty\right\}.$$
			Obviously the sequence $(001(01)^{m-2}001(01)^m)^\infty\in\widetilde{K_\beta(t_*)}$. That is $\widetilde{K_\beta(t_*)}$ contains a sequence of the smallest period 4m+2.\\
			By the definition of $S_\beta(4m+2)$ we have $S_\beta(4m+2)\geq t_*.$\\
			2) Next to prove that $S_\beta(4m+2)\leq [(001(01)^{m-2}001(01)^m)^\infty]_\beta=: t_*$. Suppose on the contrary that $S_\beta(4m+2)>t_*$.
			Then there exists $t\in (t_*,S_\beta(4m+2))$, and then $\widetilde{K_\beta(t)}$ contains at least a periodic sequence $(\omega_i)$ of the smallest period 4m+2. And we have 
			\begin{equation}
				(001(01)^{m-2}001(01)^m)^\infty\prec(t_i)\preceq \sigma^j((\omega_i))\prec(10)^\infty,\forall j\in\mathbb{N}.
			\end{equation}
			\begin{itemize}
				\item[Case 1:]
				{
					If $\omega_1=0$\\
					Since $\sigma^j((\omega_i))\prec(10)^\infty$ and $(\omega_i)$ is of the smallest period 4m+2, we can write
					$$(\omega_i)=(\square\ 1\cdots\square\ 1\ \square)^\infty,$$
					where each $\square$ only contains 0 and each $\square$ contains at least one 0.
					Then $\exists a \in \mathbb{N}$ such that $$\sigma^a((\omega_i))=(\square\ 1\cdots\square\ 1)^\infty.$$
					Since $(0101\cdots01)^\infty$ is of the smallest period 2 which is a contraction, then there is at least one $\square$ contains more than one 0. Hence $\exists b \in \mathbb{N}$ such that $$\sigma^b((\omega_i))=(00\gamma_3\gamma_4\cdots\gamma_{4m+2})^\infty.$$
					Since the inequality (10), we have $\gamma_3=1$ and $\gamma_4=0$.
					And using the inequality (10) again, we have $\gamma_5=1$ and $\gamma_6=0$. Continue this argument, we have $\gamma_{2m-1}=1$ and $\gamma_{2m}=0$.
					\begin{itemize}
						\item[Subcase 1:]{If $\gamma_{2m+1}=0.$\\
							By the inequality (3.4), we have $\gamma_{2m+2}=1$ and $\gamma_{2m+3}=0$. And then we have $\gamma_{2m+4}=1$ and $\gamma_{2m+5}=0$.
							Continue this argument, we get $\gamma_{4m+1}=0$ and then 
							\begin{align}
								\sigma^b((\omega_i))&=(001(01)^{m-2}001(01)^{m-1}0\gamma_{4m+2})^\infty\nonumber\\
								&\preceq(001(01)^{m-2}001(01)^m)^\infty \nonumber
							\end{align}
							which is a contraction to (10).}
						\item[Subcase 2:]
						{If $\gamma_{2m+1}=1,$ then we have $\gamma_{2m+2}=0.$\\
							(i) If $\gamma_{2m+3}=0$, we have $$\sigma^b((\omega_i))=(001(01)^{m-1}00\cdots \gamma_{4m+2})^\infty.$$
							By the inequality (10), we have $\gamma_{2m+4}=1$ and then $\gamma_{2m+5}=0$.\\
							\begin{itemize}
								\item[$\bullet$] 
								{
									If $\gamma_{2m+6}=0$, then we have $$\sigma^{b+2m+1}((\omega_i))=(00100\cdots\gamma_{4m+2}\gamma_{1}\cdots\gamma_{2m+1})^\infty\prec(001(01)^{m-2}001(01)^m)^\infty$$
									which is a contraction. Hence we have $\gamma_{2m+6}=1$ and then $\gamma_{2m+7}=0$.
								}
								\item [$\bullet$]
								{
									If $\gamma_{2m+8}=0$, then we have $$\sigma^{b+2m+1}((\omega_i))=(0010100\cdots\gamma_{4m+2}\gamma_{1}\cdots\gamma_{2m+1})^\infty\prec(001(01)^{m-2}001(01)^m)^\infty$$
									which is a contraction. Hence $\gamma_{2m+8}=1$ and then $\gamma_{2m+9}=0$.
								}
								
								\item[$\bullet$]
								{
									Continue this argument until  we have $\gamma_{4m-2}=1$, $\gamma_{4m-1}=0$ and $$\sigma^b((\omega_i))=(001(01)^{m-1}001(01)^{m-3}0\gamma_{4m}\gamma_{4m+1}\gamma_{4m+2})^\infty.$$
								}
								
								\item [$\bullet$]
								{
									If $\gamma_{4m}=0$, then we have $$\sigma^{b+2m+1}((\omega_i))=(001(01)^{m-3}00\cdots\gamma_{4m+2}\gamma_{1}\cdots\gamma_{2m+1})^\infty\prec(001(01)^{m-2}001(01)^m)^\infty$$
									which is a contraction. Hence we have $\gamma_{4m}=1$ and then $\gamma_{4m+1}=0$.
								}
								
								\item [$\bullet$]
								{
									If $\gamma_{4m+2}=0$, we have $$\sigma^{b+4m}((\omega_i))=(00001\cdots)^\infty\prec(001(01)^{m-2}001(01)^m)^\infty$$
									which is a contraction. Hence we have $\gamma_{4m+2}=1.$ \\
									However we have $\sigma^b=(001(01)^{m-1}001(01)^{m-1})^\infty$, which is of the smallest period 2m+1 which is a contraction.
								}
							\end{itemize}
							
							Hence there is no such (4m+2)-period sequence $(\omega_i)$ with $$\sigma^b((\omega_i))=(001(01)^{m-1}00\gamma_{2m+4}\cdots\gamma_{4m+2})^\infty.$$Therefore we have $\gamma_{2m+3}=1$ and then $\gamma_{2m+4}=0$.\\
							(ii) If $\gamma_{2m+5}=0$, we have $$\sigma^{b+2m+3}((\omega_i))=(00\gamma_{2m+6}\cdots\gamma_{4m+2}001(01)^m)^\infty\preceq(001(01)^{m-2}001(01)^m)^\infty,$$ which is a contraction to (10). Thus we have $\gamma_{2m+5}=1$ and then $\gamma_{2m+6}=0$.\\
							(iii)If $\gamma_{2m+7}=0$, we have
							\begin{align}
								\sigma^{b+2m+5}((\omega_i))&=(00\gamma_{2m+8}\cdots\gamma_{4m+2}001(01)^{m+1})^\infty\nonumber\\
								&\prec(001(01)^{m-3}001(01)^{m+1})^\infty \nonumber\\
								&\prec(001(01)^{m-2}001(01)^m)^\infty \nonumber
							\end{align}
							which is a contraction to (10). Hence we have $\gamma_{2m+7}=1$ and then $\gamma_{2m+8}=0$.\\
							Continue the same argument as $\gamma_{2m+5}=0,\gamma_{2m+7}=0$ for $\gamma_{2m+9},\gamma_{2m+11}$ and so on, we will have $\gamma_{4m-1}=1$ and corresponding $\gamma_{4m}=0$.\\
							(iv) If $\gamma_{4m+1}=0$, we have $$\sigma^{b+4m}((\omega_i))=(00\gamma_{4m+2}001(01)^{2m-2})^\infty\prec(001(01)^{m-2}001(01)^m)^\infty,$$ which is a contraction. Hence we have $\gamma_{4m+1}=1$ and then $\gamma_{4m+2}=0$. Immediately we have $$\sigma^{b+4m+1}((\omega_i))=(0001(01)^{2m-1})^\infty\prec(001(01)^{m-2}001(01)^m)^\infty$$
							which is a contraction.\\
						}
						Hence there is no such (4m+2)-period sequence $(\omega_i)$ with $\gamma_{2m+1}=1$.   			 
					\end{itemize}
				} 
				And hence there is no such (4m+2)-period sequence $(\omega_i)$ with $\omega_1=0$.
				\item[Case 2:]
				{
					If $\omega_1=1$. Since the inequality (10), we have $\omega_2$=0. Hence we have
					$$\sigma((\omega_i))=(01\ \square\ 1\cdots\square\ 1\ \square)^\infty,$$ where each $\square$ only contains 0.
					There exists $c\in \mathbb{N}$, such that $$\sigma^c((\omega_i))=(\square\ 1\ \square\ 1\cdots\square\ 1)^\infty.$$
					We transform the Case 2 into the Case 1, and corresponding we will get a contraction that there is no such (4m+2)-period sequence $(\omega_i)$ with $\omega_1=1$.
				}
			\end{itemize}
			Hence there is no such (4m+2)-sequence $(\omega_i)\in\widetilde{K_\beta(t)}$. Therefore we have $$S_\beta(4m+2)\leq t_*.$$
			Combining \ding{192} and \ding{193}, we obtain that $S_\beta(4m+2)=[(001(01)^{m-2}01(001)^m)^\infty]_\beta$.
		\end{itemize}
	\end{proof}
\end{lem}

\begin{them}\label{th:02}
	{For $\beta=\phi_2,$ we have $S_\beta(2)=0,S_\beta(4)=[(0001)^\infty]_\beta,S_\beta(6)=[(00(01)^2)^\infty]_\beta,$ and
		\begin{equation*}
			S_\beta(p)= \left\{ 
			\begin{aligned}
				&[(0(01)^m)^\infty]_\beta ,p=2m+1,m\in \mathbb{N},\\
				&[(001(01)^{m-2}001(01)^m)^\infty]_\beta ,p=4m+2,m\in \mathbb{N}_{\geq2},\\
				&[(001(01)^{m-1}001(01)^m)^\infty]_\beta ,p=4m+4,m\in \mathbb{N}_{\geq1}.\\
			\end{aligned}  
			\right.         
	\end{equation*}}
	\begin{proof}
		We have get some results based on different value of the period p. Next we will prove the last case that for $p=4m+4,m\in \mathbb{N},S_\beta(4m+4)=[(001(01)^{m-1}001(01)^m)^\infty]_\beta.$
		\begin{itemize}
			\item[(1)] {$S_\beta(4)=[(0001)^\infty]_\beta.$\\
				Firstly, we will prove that $S_\beta(4)\geq [(0001)^\infty]_\beta=: t_*$.\\
				We get that  $$\widetilde{K_\beta(t_*)}=\left\{(d_i)\in \left\{0,1\right\}^\mathbb{N} :(0001)^\infty\preceq\sigma^n((d_i))\prec(10)^\infty\right\}.$$
				Obviously the sequence $(0001)^\infty\in\widetilde{K_\beta(t_*)}$. \par 
				That is $\widetilde{K_\beta(t_*)}$ contains a sequence of the smallest period 4.\\
				By the definition of $S_\beta(4)$ we have $S_\beta(4)\geq t_*.$\\
				Next to prove that $S_\beta(4)\leq [(0001)^\infty]_\beta=: t_*$.\par Suppose on the contrary that $S_\beta(4)>t_*$.
				Then there exists $t\in (t_*,S_\beta(4))$, and then $\widetilde{K_\beta(t)}$ contains at least a periodic sequence $(\omega_i)$ of the smallest period 4. And we have 
				\begin{equation}
					(0001)^\infty\prec(t_i)\preceq \sigma^j((\omega_i))\prec(10)^\infty.
				\end{equation}
				Since $(\omega_i)$ is of the smallest period 4, there is at least one 1 among $\omega_1,\omega_2,\omega_3$ and $\omega_4$.
				And since $\sigma^j((\omega_i))\prec(10)^\infty,\forall j \in \mathbb{N}$, there is no consecutive 1 in $(\omega_i)$.
				Due to a fact that if there are more than two 1 among $\omega_1,\omega_2,\omega_3$ and $\omega_4$ there will exist consecutive 1 in $(\omega_i)$, we have 
				$$\max\# \left\{i:\omega_i=1,i\in\left\{1,2,3,4\right\}\right\}=2.$$
				If there are two 1 among $\omega_1,\omega_2,\omega_3$ and $\omega_4$, due to $\sigma^j((\omega_i))\prec(10)^\infty,\forall j \in \mathbb{N}$,\\
				we have $(\omega_i)=(1010)^\infty$ or $(\omega_i)=(0101)^\infty$ which are of the smallest period 2 which is a contraction.
				Hence there is exact one 1 among $\omega_1,\omega_2,\omega_3$ and $\omega_4$.
				Thus $\exists k\in\mathbb{N}$ such that $\sigma^k((\omega_i))=(0001)^\infty$ which is a contraction to (11).
				Therefore we have $S_\beta(4)=[(0001)^\infty]_\beta.$
			}
			\item[(2)]For $n=2m+2,m\in \mathbb{N}_{\geq1},S_\beta(4m+4)=[(001(01)^{m-1}001(01)^m)^\infty]_\beta.$\par
			\ding{192} Firstly, we will prove that $S_\beta(4m+4)\geq [(001(01)^{m-1}001(01)^m)^\infty]_\beta=: t_*$.\\
			We get that  $$\widetilde{K_\beta(t_*)}=\left\{(d_i)\in \left\{0,1\right\}^\mathbb{N} :(001(01)^{m-1}001(01)^m)^\infty\preceq\sigma^n((d_i))\prec(10)^\infty\right\}.$$
			Obviously the sequence $(001(01)^{m-1}001(01)^m)^\infty\in\widetilde{K_\beta(t_*)}$. That is $\widetilde{K_\beta(t_*)}$ contains a sequence of the smallest period 4m+4.\par
			By the definition of $S_\beta(4m+4)$ we have $S_\beta(4m+4)\geq t_*.$\par
			\ding{193} Next to prove that $S_\beta(4m+4)\leq [(001(01)^{m-1}001(01)^m)^\infty]_\beta=: t_*$.\\
			Suppose on the contrary that $S_\beta(4m+4)>t_*$.
			Then there exists $t\in (t_*,S_\beta(4m+4))$, and then $\widetilde{K_\beta(t)}$ contains at least a periodic sequence $(\omega_i)$ of the smallest period 4m+4. And we have 
			\begin{equation}
				(001(01)^{m-1}001(01)^m)^\infty\prec(t_i)\preceq \sigma^j((\omega_i))\prec(10)^\infty.
			\end{equation}
			\begin{itemize}
				\item[Case 1:]
				{
					If $\omega_1=0$.\par
					Since $\sigma^j((\omega_i))\prec(10)^\infty,\forall j \in \mathbb{N}$, there is no consecutive 1 among $\omega_1,\omega_2\cdots\omega_{4m+4}$. We can write $$(\omega_i)=(\omega_1\omega_2\cdots\omega_{4m+4})^\infty=(\square\ 1\ \square\ 1\cdots\square\ 1\ \square)^\infty,$$
					where each $\square$ only contains 0. If each $\square$ only contains one 0, we get $(\omega_i)=(01)^\infty$ which is a contraction. Hence there exists at least one $\square$ contains more than one 0. Therefore $\exists b\in \mathbb{N}$ such that $$\sigma^b((\omega_i))=(00\gamma_3\gamma_4\cdots\gamma_{4m+4})^\infty.$$
					By the inequality (12), we have $\sigma^b((\omega_i))=(001(01)^{m-1}0\gamma_{2m+3}\cdots\gamma_{4m+4})^\infty$.
					\begin{itemize}
						\item [Subcase 1:]
						{
							If $\gamma_{2m+3}=0$.\\
							Now $\sigma^b((\omega_i))=(001(01)^{m-1}00\gamma_{2m+4}\cdots\gamma_{4m+4})^\infty$. Then we have $\gamma_{2m+4}=1$ and $\gamma_{2m+5}=0$.
							\begin{itemize}
								\item [$\bullet$]
								{
									If $\gamma_{2m+6}=0$, then $$\sigma^{b+2m+1}((\omega_i))=(00100\gamma_{2m+7}\cdots\gamma_{4m+4}\gamma_1\cdots\gamma_{2m+1})^\infty\prec(001(01)^{m-1}001(01)^m)^\infty,$$which is a contraction. Hence $\gamma_{2m+6}=1$ and then $\gamma_{2m+7}=0$.
								}
								\item [$\bullet$]
								{
									If $\gamma_{2m+8}=0$, then $$\sigma^{b+2m+1}((\omega_i))=(0010100\gamma_{2m+9}\cdots\gamma_{4m+4}\gamma_1\cdots\gamma_{2m+1})^\infty\prec(001(01)^{m-1}001(01)^m)^\infty,$$ which is a contraction. Hence $\gamma_{2m+8}=1$ and then $\gamma_{2m+9}=0$.
								}\\
								Continue the same argument, we will get $\gamma_{4m+2}=1,\gamma_{4m+3}=0$ and $$\sigma^b((\omega_i))=(001(01)^{m-1}001(01)^{m-1}0\gamma_{4m+4})\preceq(001(01)^{m-1}001(01)^m)^\infty,$$
								which is a contraction. 
							\end{itemize}	
							Hence there is no such (4m+4)-period sequence $(\omega_i)$ with $\gamma_{2m+3}=0$.			
						}
						\item[Subcase 2:]
						{
							If $\gamma_{2m+3}=1$.\\
							Then we have $\gamma_{2m+4}=0.$
							\begin{itemize}
								\item[$\bullet$] 
								{
									If $\gamma_{2m+5}=0$, $\sigma^b((\omega_i))=(001(01)^{m}00\gamma_{2m+6}\cdots\gamma_{4m+4})^\infty$. Hence $$\sigma^{b+2m+3}((\omega_i))=(00\gamma_{2m+6}\cdots\gamma_{4m+4}001(01)^{m})\preceq(001(01)^{m-1}001(01)^{m})^\infty,$$which is a contraction to (12). Hence we have $\gamma_{2m+5}=1$ and then $\gamma_{2m+6}=0$.
								}
								\item [$\bullet$] 
								{
									If $\gamma_{2m+7}=0$, $\sigma^b((\omega_i))=(001(01)^{m+1}00\gamma_{2m+8}\cdots\gamma_{4m+4})^\infty$. Hence
									\begin{align}
										\sigma^{b+2m+5}((\omega_i))&=(00\gamma_{2m+8}\cdots\gamma_{4m+4}001(01)^{m+1})\nonumber\\
										&\preceq(001(01)^{m-2}001(01)^{m+1})^\infty \nonumber\\
										&\prec(001(01)^{m-1}001(01)^m)^\infty \nonumber
									\end{align} 
									which is a contraction to (12). Hence we have $\gamma_{2m+7}=1$ and then $\gamma_{2m+8}=0$.\par
								}
								Continue the same argument, we will get $\gamma_{4m+1}=1,\gamma_{4m+2}=0$ and $$\sigma^b((\omega_i))=(001(01)^{2m-1}0\gamma_{4m+3}\gamma_{4m+4})^\infty.$$
								\item[$\bullet$]
								{
									If $\gamma_{4m+3}=0$, $\sigma^b((\omega_i))=(001(01)^{2m-1}00\gamma_{4m+4})^\infty$. Hence
									\begin{align}
										\sigma^{b+4m+1}((\omega_i))&=(00\gamma_{4m+4}001(01)^{2m-1})\nonumber\\
										&\preceq(001001(01)^{2m-1})^\infty \nonumber\\
										&\prec(001(01)^{m-1}001(01)^m)^\infty\nonumber
									\end{align} 
									 which is a contraction to the inequality (12). Hence we have $\gamma_{4m+3}=1$ and then $\gamma_{4m+4}=0$.
									However at this time
									\begin{align}
										\sigma^{b+4m+3}((\omega_i))&=(0001(01)^{2m})^\infty\nonumber\\
										&\prec(001(01)^{m-1}001(01)^m)^\infty \nonumber
									\end{align} 
									which is a contraction.
								}
							\end{itemize}
							Hence there is no such (4m+4)-period sequence $(\omega_i)$ with $\gamma_{2m+3}=1$.
						}
					\end{itemize}
					Hence there is no such (4m+4)-period sequence $(\omega_i)$ with $\omega_1=0$.
				}
				\item[Case 2:]
				{
					If $\omega_1=1.$\\
					Since $\sigma^j((\omega_i))\prec(10)^\infty,\forall j\in \mathbb{N}$, we have $\omega_2=0$. Hence we have $$\sigma((\omega_i))=(01\ \square\ 1\cdots\square\ 1\ \square)^\infty.$$
					There $\exists c\in \mathbb{N}$ such that $$\sigma^c((\omega_i))=(\square\ 1\ \square\cdots\square\ 1)^\infty,$$where each $\square$ only contains 0. We transform the Case 2 into the Case 1 and correspondingly we will get a contraction. Hence there is no such (4m+4)-period sequence $(\omega_i)$ with $\omega_1=1$. 
				}
			\end{itemize}
			Hence there is no such (4m+4)-period sequence $(\omega_i)\in\widetilde{K_\beta(t)}$. \\Therefore we have $S_\beta(4m+4)\leq t_*$ and furthermore we have $S_\beta(4m+4)= t_*$.
		\end{itemize}
		Combining (1) and (2), we obtain that $S_\beta(4m+4)=[(001(01)^{m-1}01(001)^m)^\infty]_\beta$.

		To sum up, we obtain that for $\beta=\phi_2,S_\beta(2)=0,\ S_\beta(4)=(0001)^\infty_\beta,S_\beta(6)=(00(01)^2)^\infty_\beta$, and\\
		\begin{equation*}
			S_\beta(p)= \left\{ 
			\begin{aligned}
				&[(0(01)^m)^\infty]_\beta ,p=2m+1,m\in \mathbb{N},\\
				&[(001(01)^{m-2}001(01)^m)^\infty]_\beta ,p=4m+2,m\in \mathbb{N}_{\geq2},\\
				&[(001(01)^{m-1}001(01)^m)^\infty]_\beta ,p=4m+4,m\in \mathbb{N}_{\geq1}.\\
			\end{aligned}  
			\right.         
	\end{equation*}\end{proof}  
\end{them}		

\begin{cor}
	By calculation, we obtain that $S_\beta(4)=\frac{1}{\beta^4-1}\approx0.1708,S_\beta(6)=\frac{\beta^2+1}{\beta^6-1}\approx0.2135$ and
	\begin{equation*}
		S_\beta(p)= \left\{ 
		\begin{aligned}
			&\frac{1-\beta^{2m}}{(\beta^{2m+1}-1)(1-\beta^2)} ,p=2m+1,m\in \mathbb{N},\\
			&\frac{1+\beta^{2m+3}-\beta^{2m+2}-\beta^{4m+1}}{(\beta^{4m+2}-1)(1-\beta^2)} ,p=4m+2,m\in \mathbb{N}_{\geq2},\\
			&\frac{1+\beta^{2m+3}-\beta^{2m+2}-\beta^{4m+3}}{(\beta^{4m+4}-1)(1-\beta^2)} ,p=4m+4,m\in \mathbb{N}_{\geq1}.\\
		\end{aligned}  
		\right.         
	\end{equation*}
	The value of $S_\beta(p)$ is increasing to $\frac{1}{\beta^3-\beta}\approx 0.38212$ as p increasing to infinity.
\end{cor}

\begin{figure}[htb]  
	\centering  
	\begin{minipage}{.33\textwidth}  
		\centering  
		\includegraphics[width=.9\linewidth]{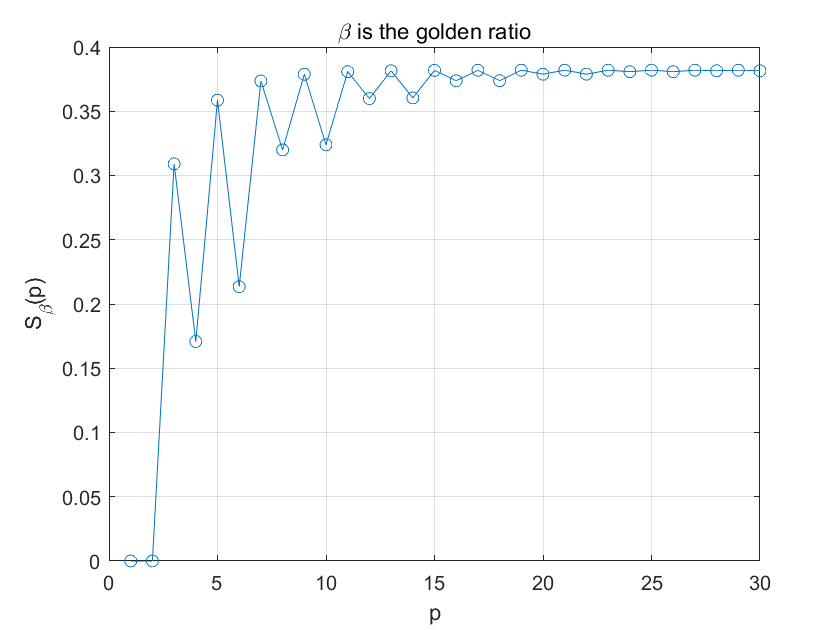} 
		\label{fig:image1}  
	\end{minipage}%
	\begin{minipage}{.33\textwidth}  
		\centering  
		\includegraphics[width=.9\linewidth]{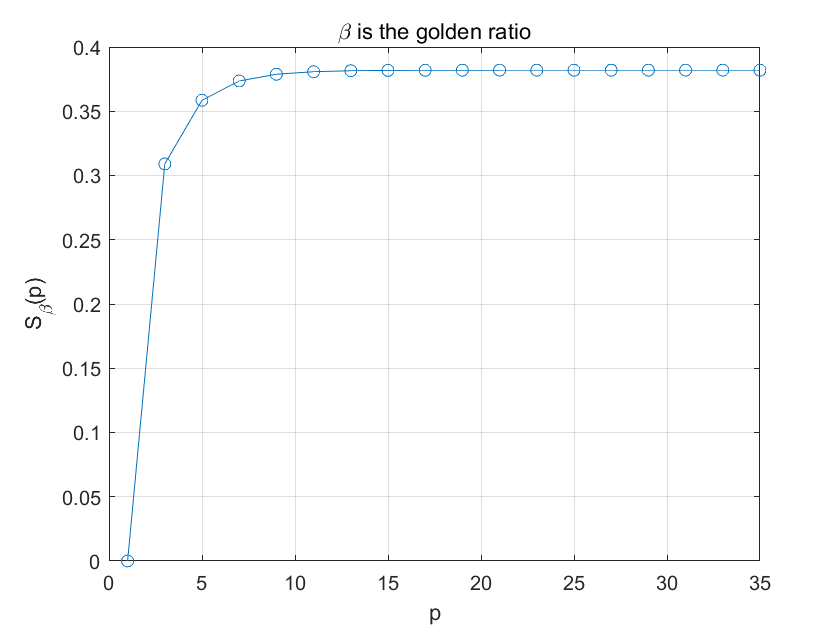} 
		\label{fig:image1}  
	\end{minipage}%
	\begin{minipage}{.33\textwidth}  
		\centering  
		\includegraphics[width=.9\linewidth]{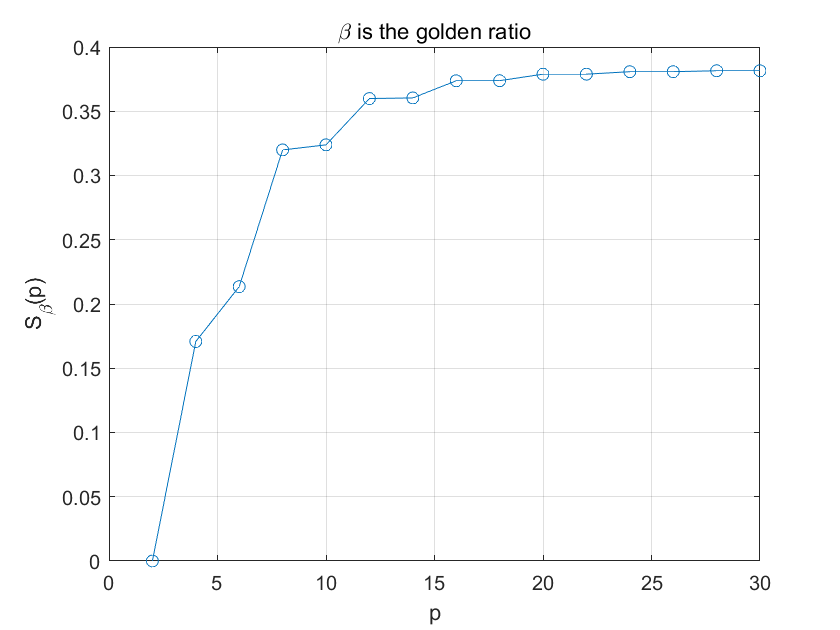}  
		\label{fig:image2}  
	\end{minipage}  
	\tiny\caption{Left:the value of $S_\beta(p)$ for $\beta=\frac{1+\sqrt{5}}{2}$ with the period p from 1 to 30.
		Middle: the value of $S_\beta(p)$ for $\beta=\frac{1+\sqrt{5}}{2}$ with the odd period p from 1 to 35. Right: the value of $S_\beta(p)$ for $\beta=\frac{1+\sqrt{5}}{2}$ with the even period p from 1 to 30.}
\end{figure} 

\section{The case for $\beta=\phi_3$}\label{sec 5}
{
	In the following section, we will consider the more complex situation that replacing the golden ratio with the tribonacci number which is the unique root in $(1,2)$ of $x^3=x^2+x+1$ which is approx to 1.8393. Importantly, if $\beta$ is a multinacci number of degree m then its quasi-greedy expansion of 1 is of the form $(1^{m-1}0)^\infty$. Corresponding to the property, the quasi-greedy expansion of 1 is $\delta(\beta)=(1^20)^\infty$ where $\beta$ is the tribonacci number which is the multinacci number of degree 3.
}

\begin{lem} 
	{For $\beta=\phi_3,\ \delta(\beta)=(1^20)^\infty$, we have $S_\beta(1)=0,S_\beta(3)=[(001)^\infty]_\beta,S_\beta(4)=[(0011)^\infty]_\beta$ and $S_\beta(6)=[(001101)^\infty]_\beta$.
	}
	\begin{proof}
		\begin{itemize}
			\item [(1)]
			{
				For $p=1$, $S_\beta(1)=0$.
				\begin{itemize}
					\item [\ding{192}]
					{
						By the definition of $S_\beta(1)$, we can get $S_\beta(1)\geq 0$ trivially.
					}
					\item [\ding{193}]
					{
						Next to prove that $S_\beta(1)\leq 0$.\\
						Suppose on the contrary that $S_\beta(1)>0$, then there exists $t\in (0,S_\beta(1))$, and then $\widetilde{K_\beta(t)}$ contains at least a periodic sequence $(\omega_i)$ of the smallest period 1. And we have 
						$$(0)^\infty\prec(t_i)\preceq \sigma^j((\omega_i))\prec(1^20)^\infty.$$
						There are only $(0)^\infty$ and $(1)^\infty$ two periodic sequences of the smallest period 1. For $(1)^\infty$ we have $(1)^\infty\succ(1^20)^\infty$ and for $(0)^\infty$ we get a contraction to the inequality above. Hence there is no such 1-period sequence $(\omega_i)\in \widetilde{K_\beta(t)}$, we get a contraction.
						Therefore we have $S_\beta(1)\leq 0$.  		 		
					}
					
				\end{itemize}
				Combined \ding{192} and \ding{193} we have $S_\beta(1)= 0.$	
			}
			\item[(2)] 
			{
				For $p=3$, $S_\beta(3)=[(001)^\infty]_\beta$.
				\begin{itemize}
					\item [\ding{192}]
					{
						Firstly, we will prove that $S_\beta(3)\geq [(001)^\infty]_\beta=: t_*$.\\
						We get that  $$\widetilde{K_\beta(t_*)}=\left\{(d_i)\in \left\{0,1\right\}^\mathbb{N} :(001)^\infty\preceq\sigma^n((d_i))\prec(1^20)^\infty\right\}.$$
						Obviously the sequence $(001)^\infty\in\widetilde{K_\beta(t_*)}$. That is $\widetilde{K_\beta(t_*)}$ contains a sequence of the smallest period 3.\\
						By the definition of $S_\beta(3)$ ,we have $S_\beta(3)\geq t_*.$
					}
					\item [\ding{193}]
					{
						Next to prove that $S_\beta(3)\leq [(001)^\infty]_\beta=: t_*$. \\
						Suppose on the contrary that $S_\beta(3)>[(001)^\infty]_\beta$, then there exists $t\in (t_*,S_\beta(3))$, and then $\widetilde{K_\beta(t)}$ contains at least a periodic sequence $(\omega_i)$ of the smallest period 3. And we have 
						\begin{equation}
							(001)^\infty\prec(t_i)\preceq \sigma^j((\omega_i))\prec(1^20)^\infty,\forall j \in \mathbb{N}.
						\end{equation}
						Since $(\omega_i)$ is of the smallest period 3, there is at least one 1 among $\omega_1,\omega_2$ and $\omega_3$.
						And since $\sigma^j((\omega_i))\prec(1^20)^\infty,\forall j \in \mathbb{N}, (\omega_i)\neq (111)^\infty$ which implies that there exists 0 among $\omega_1,\omega_2$ and $\omega_3$.
						Therefore $\exists a \in \mathbb{N}$ such that $\sigma^a((\omega_i))=(0x_2x_3)^\infty.$
						\begin{itemize}
							\item [$\bullet$]
							{
								If $x_2=0, \sigma^a((\omega_i))=(00x_3)^\infty\preceq(001)^\infty$ which is a contraction to (13). Hence $x_2=1$.
							}
							\item [$\bullet$]
							{
								If $x_3=1, \sigma^a((\omega_i))=(011)^\infty$ and $\sigma^{a+1}((\omega_i))=(1^20)^\infty$ which is a contraction to $\sigma^j((\omega_i))\prec(1^20)^\infty,\forall j \in \mathbb{N}$. Hence we have $x_3=0$.
								However at this time $\sigma^{a+2}=(001)^\infty$ which is also a contraction to (13).
							}
						\end{itemize}
						Hence there is no such 3-period sequence $(\omega_i)\in \widetilde{K_\beta(t)}$.
						Therefore we have $$S_\beta(3)\leq [(001)^\infty]_\beta.$$		
					}
					
				\end{itemize}
				Combined \ding{192} and \ding{193} we have $S_\beta(3)= [(001)^\infty]_\beta.$
			}
			\item [(3)]
			{
				For p=4, $S_\beta(4)=[(0011)^\infty]_\beta$.
				\begin{itemize}
					\item [\ding{192}]
					{
						Firstly, we will prove that $S_\beta(4)\geq [(0011)^\infty]_\beta=: t_*$.\\
						We get that  $$\widetilde{K_\beta(t_*)}=\left\{(d_i)\in \left\{0,1\right\}^\mathbb{N} :(0011)^\infty\preceq\sigma^n((d_i))\prec(1^20)^\infty\right\}.$$
						Obviously the sequence $(0011)^\infty\in\widetilde{K_\beta(t_*)}$. That is $\widetilde{K_\beta(t_*)}$ contains a sequence of the smallest period 4.\\
						By the definition pf $S_\beta(4)$ , we have $S_\beta(4)\geq t_*.$
					}
					\item [\ding{193}]
					{
						Next to prove that $S_\beta(4)\leq [(0011)^\infty]_\beta=: t_*$. Suppose on the contrary that $S_\beta(4)>[(0011)^\infty]_\beta$.
						Then there exists $t\in (t_*,S_\beta(4))$, and then $\widetilde{K_\beta(t)}$ contains at least a periodic sequence $(\omega_i)$ of the smallest period 4. And we have
						\begin{equation}
							(0011)^\infty\prec(t_i)\preceq \sigma^j((\omega_i))\prec(1^20)^\infty.
						\end{equation}           
						Since $(\omega_i)$ is of the smallest period 4, there is at least one 1 among $\omega_1,\omega_2,\omega_3$ and $\omega_4$.
						And since $\sigma^j((\omega_i))\prec(1^20)^\infty,\forall j \in \mathbb{N}$, there is no consecutive three 1 in $(\omega_i)$  which implies that there exists 0 among $\omega_1,\omega_2,\omega_3$ and $\omega_4$.
						Therefore $\exists a \in \mathbb{N}$ such that $\sigma^a((\omega_i))=(0x_2x_3x_4)^\infty.$
						\begin{itemize}
							\item [$\bullet$]
							{
								If $x_2=0, \sigma^a((\omega_i))=(00x_3x_4)^\infty\preceq(0011)^\infty$ which is a contraction to (14). Hence $x_2=1$.
							}
							\item [$\bullet$]
							{
								If $x_3=1$, then $x_4=0$ $\sigma^{a}((\omega_i))=(0110)^\infty$ and $\sigma^{a+3}((\omega_i))=(0011)^\infty$ which is a contraction to (14). Hence we have $x_3=0$.
							}
							\item[$\bullet$]
							{
								If $x_4=1, \sigma^a((\omega_i))=(0101)^\infty$ which is of the smallest period 2 which is a contraction. Hence $x_4=0$. However $\sigma^a((\omega_i))=(0100)^\infty$ which is contradict to $\sigma^{a+2}((\omega_i))\succ(0011)^\infty$. 
							}
						\end{itemize}
						Hence there is no such 4-period sequence $(\omega_i)\in \widetilde{K_\beta(t)}$.
						Therefore we have $$S_\beta(4)\leq [(0011)^\infty]_\beta.$$
					}
				\end{itemize}
				Combined \ding{192} and \ding{193} we have $S_\beta(4)= [(0011)^\infty]_\beta.$
			}
			\item [(4)]
			{
				For p=6, $S_\beta(6)=[(001101)^\infty]_\beta$.
				\begin{itemize}
					\item [\ding{192}]
					{
						Firstly, we will prove that $S_\beta(6)\geq [(001101)^\infty]_\beta=: t_*$.\\
						We get that  $$\widetilde{K_\beta(t_*)}=\left\{(d_i)\in \left\{0,1\right\}^\mathbb{N} :(001101)^\infty\preceq\sigma^n((d_i))\prec(1^20)^\infty\right\}.$$
						Obviously the sequence $(001101)^\infty\in\widetilde{K_\beta(t_*)}$. That is $\widetilde{K_\beta(t_*)}$ contains a sequence of the smallest period 6.\\
						By the definition of $S_\beta(6)$ , we have $S_\beta(6)\geq t_*.$
					}
					\item [\ding{193}]
					{
						Next to prove that $S_\beta(6)\leq [(001101)^\infty]_\beta=: t_*$. \\
						Suppose on the contrary that $S_\beta(6)>t_*$.
						Then there exists $t\in (t_*,S_\beta(6))$, and then $\widetilde{K_\beta(t)}$ contains at least a periodic sequence $(\omega_i)$ of the smallest period 6. And we have
						\begin{equation}
							(001101)^\infty\prec(t_i)\preceq \sigma^j((\omega_i))\prec(1^20)^\infty, \forall j \in \mathbb{N}.
						\end{equation} 
						Since $(\omega_i)$ is of the smallest period 6, there is at least one 1 among $\omega_1,\omega_2,\omega_3\cdots\omega_6$.
						And since $\sigma^j((\omega_i))\prec(1^20)^\infty,\forall j \in \mathbb{N}$, there is no consecutive three 1 in $(\omega_i)$  which implies that there exists 0 among $\omega_1,\omega_2,\omega_3\cdots\omega_6$.
						Therefore $\exists a \in \mathbb{N}$ such that $\sigma^a((\omega_i))=(0x_2x_3x_4x_5x_6)^\infty.$\par
						{
							If $x_2=0$, we have $x_3=x_4=1$ which is due to the inequality (15) and then $x_5=0$. $\sigma^a((\omega_i))=(00110x_6)^\infty\preceq(001101)^\infty$ which is a contraction to (15). \par 
							Hence $x_2=1$.
						}
						\begin{itemize}
							
							\item [(i)]
							{
								If $x_3=1$, then $x_4=0$ and $\sigma^{a}((\omega_i))=(0110x_5x_6)^\infty$.
								\begin{itemize}
									\item [$\bullet$]
									{
										If $x_5=0, \sigma^{a+3}((\omega_i))=(00x_6011)^\infty\prec(001101)^\infty$ which is a contraction to (15). Hence we have $x_5=1$.
									}
									\item[$\bullet$]
									{
										If $x_6=1,\sigma^a((\omega_i))=(011011)^\infty$ which is of the smallest period 3 which is a contradiction. Hence $x_6=0$. But $\sigma^a((\omega_i))=(011010)$, we have $$\sigma^{a+5}((\omega_i))=(001101)^\infty$$ which is a contraction to (15).
									}
								\end{itemize}
								
							}
							\item[(ii)]
							{
								If $x_3=0$, then $\sigma^a((\omega_i))=(010x_4x_5x_6)^\infty$.
								\begin{itemize}
									\item [$\bullet$]
									{
										If $x_4=0$, we have $\sigma^{a+2}((\omega_i))=(00x_5x_601)^\infty\preceq(001101)^\infty$ which is a contraction.
										Hence $x_4=1$.
									}
									\item [$\bullet$]
									{
										If $x_5=1$, then $x_6=0$. We have $$\sigma^{a+5}((\omega_i))=(001011)^\infty\prec(001101)^\infty$$ which is a contraction to (15). Hence $x_5=0$.
									}
									\item[$\bullet$]
									{
										If $x_6=1, \sigma^a((\omega_i))=(010101)^\infty$ which of the smallest period 2 which is a contraction. Hence $x_6=0$. But $\sigma^{a+4}((\omega_i))=(000101)^\infty\prec(001101)^\infty$ which is also a contraction. 
									}
								\end{itemize} 
							}
							
						\end{itemize}
						Hence there is no such 6-period sequence $(\omega_i)\in \widetilde{K_\beta(t)}$ which is a contraction. Therefore we have $S_\beta(6)\leq [(001101)^\infty]_\beta$.
					}
				\end{itemize}
				Combined \ding{192} and \ding{193} we have $S_\beta(6)= [(001101)^\infty]_\beta.$
			}
		\end{itemize}
	\end{proof}
\end{lem}

\begin{lem}
	{For $\beta=\phi_3,\ \delta(\beta)=(1^20)^\infty,$ we have $$\forall m\in \mathbb{N},S_\beta(3m+2)=[(01(011)^m)^\infty]_\beta.$$
	}
	\begin{proof}
		\begin{itemize}
			\item [\ding{192}]
			{
				Firstly, we will prove that $S_\beta(3m+2)\geq [(01(011)^m)^\infty]_\beta=: t_*$.\\
				We get that  $$\widetilde{K_\beta(t_*)}=\left\{(d_i)\in \left\{0,1\right\}^\mathbb{N} :(01(011)^m)^\infty\preceq\sigma^n((d_i))\prec(1^20)^\infty\right\}.$$
				Obviously the sequence $(01(011)^m)^\infty\in\widetilde{K_\beta(t_*)}$. That is $\widetilde{K_\beta(t_*)}$ contains a sequence of the smallest period 3m+2.\\
				By the definition of $S_\beta(3m+2)$ , we have $S_\beta(3m+2)\geq t_*.$
			}
			\item [\ding{193}]
			{
				Next to prove that $S_\beta(3m+2)\leq [(01(011)^m)^\infty]_\beta=: t_*$.\\
				Suppose on the contrary that $S_\beta(3m+2)>t_*$.
				Then there exists $t\in (t_*,S_\beta(3m+2))$, and then $\widetilde{K_\beta(t)}$ contains at least a periodic sequence $(\omega_i)$ of the smallest period 3m+2. And we have 
				\begin{equation}
					(01(011)^m)^\infty\prec(t_i)\preceq \sigma^j((\omega_i))\prec(1^20)^\infty, \forall j \in \mathbb{N}.
				\end{equation}
				Since $(\omega_i)$ is of the smallest period 3m+2, there is at least one 1 among $\omega_1,\omega_2,\omega_3\cdots\omega_{3m+2}$.
				And since $\sigma^j((\omega_i))\prec(1^20)^\infty,\forall j \in \mathbb{N}$, there is no consecutive three 1 in $(\omega_i)$  which implies that there exists 0 among $\omega_1,\omega_2,\omega_3\cdots\omega_{3m+2}$.
				Therefore $\exists a \in \mathbb{N}$ such that $$\sigma^a((\omega_i))=(0x_2x_3\cdots x_{3m+2})^\infty.$$
				Since $\sigma^j((\omega_i))\succ(01(011)^m)^\infty,\forall j\in \mathbb{N}$, we have $x_2=1$.
				\begin{itemize}
					\item[Case 1:]
					{
						If $x_3=0,\sigma^a((\omega_i))=(010x_4\cdots x_{3m+2})\preceq(01(011)^m)^\infty$ which is a contraction to (16).
						Hence there is no such (3m+2)-period sequence $(\omega_i)\in \widetilde{K_\beta(t)}$ with $x_3=0$.
					}
					\item [Case 2:]
					{
						If $x_3=1$, then $x_4=0$. Since $\sigma^j((\omega_i))\succ(01(011)^m)^\infty$, there is no consecutive 0 in $(\omega_i)$. Hence we have $x_5=1$.
						\begin{itemize}
							\item[$\bullet$]
							{
								If $x_6=0, \sigma^{a+3}((\omega_i))=(010x_7\cdots x_{3m+2}011)^\infty\preceq(01(011)^m)^\infty$, which is a contraction to (16). Hence we have $x_6=1$ then $x_7=0$ and $x_8=1$.
							} 
							\item[$\bullet$]
							{
								If $x_9=0,\sigma^{a+6}((\omega_i))=(010x_{10}\cdots x_{3m+2}011011)^\infty\preceq(01(011)^m)^\infty$, which is a contraction to (16). Hence we have $x_9=1$ then $x_{10}=0$ and $x_{11}=1$. 
							} \\
							Continue the same argument, we will get that $x_{3m-3}=1$ then $x_{3m-2}=0$ and $x_{3m-1}=1$.
							\item[$\bullet$]
							{
								If $x_{3m}=0,\sigma^{a+3m-3}((\omega_i))=(010x_{3m+1}x_{3m+2}(011)^{m-1})^\infty\preceq(01(011)^m)^\infty$, which is a contraction to (16). Hence we have $x_{3m}=1$ then $x_{3m+1}=0$ and $x_{3m+2}=1$. However at this time $\sigma^{a+3m}((\omega_i))=(01(011)^m)$which is a contraction to (16). 
							}
						\end{itemize}
						Hence there is no such (3m+2)-period sequence $(\omega_i)\in \widetilde{K_\beta(t)}$ with $x_3=1$ or $x_3=0$ which is a contraction.
					}
				\end{itemize}
				Hence we have $S_\beta(3m+2)\leq[(01(011)^m)^\infty]_\beta$.
			}
			
		\end{itemize}
		Combined \ding{192} and \ding{193} we have $S_\beta(3m+2)= [(01(011)^m)^\infty]_\beta.$
	\end{proof}
\end{lem}

\begin{cor}
	By calculation, we obtain that $$S_\beta(3m+2)=\frac{1+\beta-\beta^{3m+1}-\beta^{3m+3}}{(1-\beta^3)(\beta^{3m+2}-1)},\forall m\in\mathbb{N}.$$
	The value of $S_\beta(p)$ is increasing to $\frac{\beta^{2}+1}{\beta^4-\beta}\approx 0.45626$ as $p$ increasing to infinity.
	
	\begin{figure}[htb]  
		\centering  
		\begin{minipage}{.5\textwidth}  
			\centering  
			\includegraphics[width=.8\linewidth]{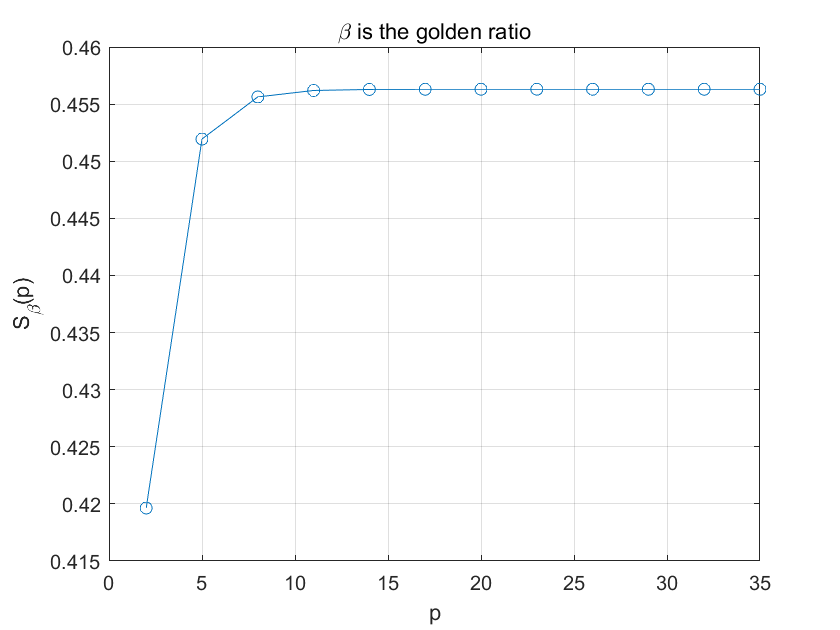}  
			\label{fig:image2}  
		\end{minipage}
		\tiny\caption{The value of $S_\beta(p)$ for $\beta\approx1.8393$ with the period $p=3m+2,\forall m\in\mathbb{N}$ from 1 to 35.}
	\end{figure} 
\end{cor}

\begin{lem}
	{For $\beta=\phi_3,\ \delta(\beta)=(1^20)^\infty,$ $\forall m\in \mathbb{N}_{\geq2}$, we have
		\begin{equation*}
			S_\beta(3m+1)= \left\{ 
			\begin{aligned}
				&[(01(011)^{z-1}01(011)^z)^\infty]_\beta ,\ m=2z,z\in \mathbb{N}_{\geq1},\\
				&[(01(011)^{z-1}01(011)^{z+1})^\infty]_\beta,\ m=2z+1,z\in \mathbb{N}_{\geq1}.\\
			\end{aligned}  
			\right.        
		\end{equation*}
	}
	\begin{proof}
		
		{(1)For $p=3m+1,m=2z,z\in \mathbb{N}_{\geq1},S_\beta(p)=[(01(011)^{z-1}01(011)^z)^\infty]_\beta$.\\
			{1)Firstly, we will prove that $$S_\beta(6z+1)\geq [(01(011)^{z-1}01(011)^z)^\infty]_\beta=: t_*,\forall z \in \mathbb{N}_{\geq1}.$$\\
				We get that  $$\widetilde{K_\beta(t_*)}=\left\{(d_i)\in \left\{0,1\right\}^\mathbb{N} :(01(011)^{z-1}01(011)^z)^\infty\preceq\sigma^n((d_i))\prec(1^20)^\infty\right\}.$$
				Obviously the sequence $(01(011)^{z-1}01(011)^z)^\infty\in\widetilde{K_\beta(t_*)}$. \par 
				That is $\widetilde{K_\beta(t_*)}$ contains a sequence of the smallest period 6z+1.\\
				By the definition of $S_\beta(6z+1)$, we have $S_\beta(6z+1)\geq t_*.$\par
			}
			{2)Next to prove that $S_\beta(6z+1)\leq [(01(011)^{z-1}01(011)^z)^\infty]_\beta=: t_*$. \\
				Suppose on the contrary that $S_\beta(6z+1)>t_*$, then there exists $t\in (t_*,S_\beta(6z+1))$, and then $\widetilde{K_\beta(t)}$ contains at least a periodic sequence $(\omega_i)$ of the smallest period 6z+1. And we have 
				\begin{equation}
					(01(011)^{z-1}01(011)^z)^\infty\prec(t_i)\preceq \sigma^j((\omega_i))\prec(1^20)^\infty.
				\end{equation}
				Since $(\omega_i)$ is of the smallest period 6z+1, there is at least one 1 among $\omega_1,\omega_2,\omega_3\cdots\omega_{6z+1}$.
				Since $\sigma^j((\omega_i))\succ(01(011)^{z-1}01(011)^z)^\infty,\forall j \in \mathbb{N}$, there is no consecutive 0 in $(\omega_i)$.
				And since $\sigma^j((\omega_i))\prec(1^20)^\infty,\forall j \in \mathbb{N}$, there is no consecutive three 1 in $(\omega_i)$  which implies that there exists 0 among $\omega_1,\omega_2,\omega_3\cdots\omega_{6z+1}$.
				Therefore $\exists a \in \mathbb{N}$ such that $$\sigma^a((\omega_i))=(0x_2x_3\cdots x_{6z+1})^\infty.$$
				Since $\sigma^j((\omega_i))\succ(01(011)^{z-1}01(011)^z)^\infty,\forall j\in \mathbb{N}$, we have $x_2=1$.
				\begin{itemize}
					\item [\ding{192}]
					{
						If $x_3=0$, since
						 $$\sigma^j((\omega_i))\succ(01(011)^{z-1}01(011)^z)^\infty,\forall j\in \mathbb{N}$$
						 we have $\sigma^a((\omega_i))=(01(011)^{z-1}01x_{3z+2}\cdots x_{6z+1})^\infty$.
						\begin{itemize}
							\item[$\bullet$]
							{
								If $x_{3z+2}=0$, then $$\sigma^a((\omega_i))=(01(011)^{z-1}010\cdots x_{6z+1})^\infty\preceq(01(011)^{z-1}01(011)^{z})^\infty$$ which is a contraction to (17). Hence we get $x_{3z+2}=1$ then $x_{3z+3}=0$ and $x_{3z+4}=1$.
							}
							\item[$\bullet$]
							{
								If $x_{3z+5}=0$, then $$\sigma^a((\omega_i))=(01(011)^{z}010\cdots x_{6z+1})^\infty.$$We have
								$$\sigma^{a+3z+2}((\omega_i))=(010x_{3z+6}\cdots x_{6z+1}01(011)^{z})^\infty\preceq(01(011)^{z-1}01(011)^{z})^\infty$$ which is a contraction to (17). Hence we get $x_{3z+5}=1$ then $x_{3z+6}=0$ and $x_{3z+7}=1$.
							}
							\item[$\bullet$]
							{
								If $x_{3z+8}=0$, then $$\sigma^a((\omega_i))=(01(011)^{z+1}010\cdots x_{6z+1})^\infty.$$We have
								$$\sigma^{a+3z+5}((\omega_i))=(010x_{3z+9}\cdots x_{6z+1}01(011)^{z+1})^\infty\prec(01(011)^{z-1}01(011)^{z})^\infty$$ which is a contraction to (17). Thus we get $x_{3z+8}=1$ then $x_{3z+9}=0$ and $x_{3z+10}=1.$
							}\\
							\item[$\bullet$]
							{Continue the same argument as $x_{3z+2}=0,x_{3z+5}=0$ and $x_{3z+8}=0$ until $x_{6z-1}=0$. We will get a contraction which implies $x_{6z-1}=1$ then $x_{6z}=0,x_{6z+1}=1$ and $$\sigma^a((\omega_i))=(01(011)^{2z-1}01)^\infty.$$
								However $\sigma^{a+6z-1}((\omega_i))=(0101(011)^{z-1})^\infty\prec(01(011)^{z-1}01(011)^z)^\infty$, which is a contraction to (17).
							}
						\end{itemize} 
						Hence there is no such (6z+1)-period sequence $(\omega_i)\in \widetilde{K_\beta(t)}$ with $x_3=0$.
					}
					\item [\ding{193}]
					{
						If $x_3=1$, then $x_4=0,x_5=1$ and 
						$$\sigma^a((\omega_i))=(01101x_6\cdots x_{6z+1})^\infty.$$
						\begin{itemize}
							\item [\textrm{(i)}]
							{
								If $x_6=0, \sigma^{a+3}((\omega_i))=:(y_i)=(010y_4\cdots y_{6z-2}011)^\infty$\\
								since $$\sigma^j((\omega_i))\succ(01(011)^{z-1}01(011)^z)^\infty,\forall j\in \mathbb{N},$$ we have $(y_i)=(01(011)^{z-1}01y_{3z+2}\cdots y_{6z-2}011)^\infty$.
								\begin{itemize}
									\item[$\bullet$]
									{
										If $y_{3z+2}=0$, then $$(y_i)=(01(011)^{z-1}010\cdots y_{6z-2}011)^\infty\preceq(01(011)^{z-1}01(011)^{z})^\infty$$ which is a contraction to (17). \\Hence we get $y_{3z+2}=1$ then $y_{3z+3}=0$ and $y_{3z+4}=1$.
									}
									\item[$\bullet$]
									{
										If $y_{3z+5}=0$, then $$(y_i)=(01(011)^{z}010\cdots y_{6z-2}011)^\infty.$$We have
										$$\sigma^{3z+2}((y_i))=(010y_{3z+6}\cdots y_{6z-2}01101(011)^{z})^\infty\preceq(01(011)^{z-1}01(011)^{z})^\infty$$ which is a contraction to (17). \\Hence we get $y_{3z+5}=1$ then $y_{3z+6}=0$ and $y_{3z+7}=1$.
									}
									\item[$\bullet$]
									{
										If $y_{3z+8}=0$, then $$(y_i)=(01(011)^{z+1}010\cdots y_{6z-2}011)^\infty.$$We have
										$$\sigma^{3z+5}((y_i))=(010y_{3z+9}\cdots y_{6z-2}01101(011)^{z+1})^\infty\prec(01(011)^{z-1}01(011)^{z})^\infty$$ which is a contraction. \\Thus we get $y_{3z+8}=1$ then $y_{3z+9}=0$ and $y_{3z+10}=1.$
									}\\
									\item[$\bullet$]
									{Continue the same argument as $y_{3z+2}=0,y_{3z+5}=0$ and $y_{3z+8}=0$ until $y_{6z-4}=0$. We will get a contraction which implies $y_{6z-4}=1$ then $y_{6z-3}=0,y_{6z-2}=1$ and $$(y_i)=(01(011)^{2z-2}01011)^\infty.$$
										However $\sigma^{6z-4}((y_i))=(0101101(011)^{2z-2})^\infty\prec(01(011)^{z-1}01(011)^z)^\infty$, which is a contraction to (17).
									}
								\end{itemize} 
								Hence there is no such (6z+1)-period sequence $(\omega_i)\in \widetilde{K_\beta(t)}$ with $y_3=0$. Hence we have $x_6=y_3=1$ then $x_7=0$ and $x_8=1$.
							}
							\item [\textrm{(ii)}]
							{
								If $x_9=0, \sigma^{a+6}((\omega_i))=:(y_i)=(010y_4\cdots y_{6z-5}011011)^\infty$.\\
								Since $\sigma^j((\omega_i))\succ(01(011)^{z-1}01(011)^z)^\infty,\forall j\in \mathbb{N}$,we have $$(y_i)=(01(011)^{z-1}01y_{3z+2}\cdots y_{6z-5}011011)^\infty.$$
								\begin{itemize}
									\item[$\bullet$]
									{
										If $y_{3z+2}=0$, then $$(y_i)=(01(011)^{z-1}010\cdots y_{6z-5}011011)^\infty\preceq(01(011)^{z-1}01(011)^{z})^\infty$$ which is a contraction to (17). \\Hence we get $y_{3z+2}=1$ then $y_{3z+3}=0$ and $y_{3z+4}=1$.
									}
									\item[$\bullet$]
									{
										If $y_{3z+5}=0$, then $$(y_i)=(01(011)^{z}010\cdots y_{6z-5}011011)^\infty.$$We have
										$$\sigma^{3z+2}((y_i))=(010y_{3z+6}\cdots y_{6z-5}01101101(011)^{z})^\infty\preceq(01(011)^{z-1}01(011)^{z})^\infty$$ which is a contraction to(17). \\Hence we get $y_{3z+5}=1$ then $y_{3z+6}=0$ and $y_{3z+7}=1$.
									}
									\item[$\bullet$]
									{
										If $y_{3z+8}=0$, then $$(y_i)=(01(011)^{z+1}010\cdots y_{6z-5}011011)^\infty.$$We have
										$$\sigma^{3z+5}((y_i))=(010y_{3z+9}\cdots y_{6z-5}01101101(011)^{z+1})^\infty\prec(01(011)^{z-1}01(011)^{z})^\infty$$ which is a contraction to (17). \\Thus we get $y_{3z+8}=1$ then $y_{3z+9}=0$ and $y_{3z+10}=1.$
									}
									\item[$\bullet$]
									{Continue the same argument as $y_{3z+2}=0,y_{3z+5}=0$ and $y_{3z+8}=0$ until $y_{6z-7}=0$. We will get a contraction which implies $y_{6z-7}=1$ then $y_{6z-6}=0,y_{6z-5}=1$ and $$(y_i)=(01(011)^{2z-3}01011011)^\infty.$$
										However $\sigma^{6z-7}((y_i))=(0101101101(011)^{2z-3})^\infty\prec(01(011)^{z-1}01(011)^z)^\infty$, which is a contraction to (17).
									}
								\end{itemize} 
								Hence there is no such (6z+1)-period sequence $(\omega_i)\in \widetilde{K_\beta(t)}$ with $x_9=0$. Hence $x_9=1$ then $x_{10}=0$ and $x_{11}=1$.
							}\\
							Continue the same argument as $x_3=0,x_6=0$ and $x_9=0$ until $x_{6z-3}=0$, we will get a contraction which implies that $x_{6z-3}=1$ then $x_{6z-2}=0,x_{6z-1}=1$ and
							$$\sigma^a((\omega_i))=((011)^{2m-1}01x_{6z}x_{6z+1})^\infty.$$
							\item [\textrm{(iii)}]
							{
								If $x_{6z}=1$, then $x_{6z+1}=0.$\\
								At this time we have $$\sigma^{a+6z}((\omega_i))=(0(011)^{2z})\prec(01(011)^{z-1}01(011)^z)^\infty$$ which is a contraction to (17). Hence $x_{6z}=0$. However  $$\sigma^{a+6z-3}((\omega_i))=(010x_{6z+1}(011)^{2z-1})\preceq(0101(011)^{2z-1})\prec(01(011)^{z-1}01(011)^z)^\infty$$ which is also a contraction to (17).
							}
						\end{itemize}
						
					}
					Hence there is no such (6z+1)-sequence $(\omega_i)\in \widetilde{K_\beta(t)}$ which is a contraction. Hence $\forall z\in \mathbb{N}_{\geq1},S_\beta(6z+1)\leq[(01(011)^{z-1}01(011)^z)^\infty]_\beta$.  
				\end{itemize}
				Combined 1) with 2), we have $$\forall z\in \mathbb{N}_{\geq1},S_\beta(6z+1)=[(01(011)^{z-1}01(011)^z)^\infty]_\beta.$$
			}
		}
		{(2)For $p=3m+1,m=2z+1,z\in \mathbb{N}_{\geq1},S_(p)=[(01(011)^{z-1}01(011)^{z+1})^\infty]_\beta$.\\
			{1)Firstly, we will prove that $S_\beta(6z+4)\geq [(01(011)^{z-1}01(011)^{z+1})^\infty]_\beta=: t_*$.\\
				We get that  $$\widetilde{K_\beta(t_*)}=\left\{(d_i)\in \left\{0,1\right\}^\mathbb{N} :(01(011)^{z-1}01(011)^{z+1})^\infty\preceq\sigma^n((d_i))\prec(1^20)^\infty\right\}.$$
				Obviously the sequence $(01(011)^{z-1}01(011)^{z+1})^\infty\in\widetilde{K_\beta(t_*)}$. \\
				That is $\widetilde{K_\beta(t_*)}$ contains a sequence of the smallest period 6z+4.\\
				By the definition of $S_\beta(6z+4)$ , we have $S_\beta(6z+4)\geq t_*.$\\
			}
			{2)Next to prove that $S_\beta(6z+4)\leq [(01(011)^{z-1}01(011)^{z+1})^\infty]_\beta=: t_*$. \\
				Suppose on the contrary that $S_\beta(6z+4)>t_*$, then there exists $t\in (t_*,S_\beta(6z+4))$, and then $\widetilde{K_\beta(t)}$ contains at least a periodic sequence $(\omega_i)$ of the smallest period 6z+4. And we have 
				\begin{equation}
					(01(011)^{z-1}01(011)^{z+1})^\infty\prec(t_i)\preceq \sigma^j((\omega_i))\prec(1^20)^\infty.
				\end{equation}
				Since $(\omega_i)$ is of the smallest period 6z+4, there is at least one 1 among $\omega_1,\omega_2,\omega_3\cdots\omega_{6z+4}$.
				And since $\sigma^j((\omega_i))\prec(1^20)^\infty,\forall j \in \mathbb{N}$, there is no consecutive three 1 in $(\omega_i)$  which implies that there exists 0 among $\omega_1,\omega_2,\omega_3\cdots\omega_{6z+4}$.
				Therefore $\exists a \in \mathbb{N}$ such that $$\sigma^a((\omega_i))=(0x_2x_3\cdots x_{6z+4})^\infty.$$
				Since $\sigma^j((\omega_i))\succ(01(011)^{z-1}01(011)^{z+1})^\infty,\forall j\in \mathbb{N}$, we have $x_2=1$.
				\begin{itemize}
					\item [\ding{192}]
					{
						If $x_3=0$, \\
						since $\sigma^j((\omega_i))\succ(01(011)^{z-1}01(011)^{z+1})^\infty,\forall j\in \mathbb{N}$,\\
						we have $\sigma^a((\omega_i))=(01(011)^{z-1}01x_{3z+2}\cdots x_{6z+4})^\infty$.
						\begin{itemize}
							\item[$\bullet$]
							{
								If $x_{3z+2}=0$, then $$\sigma^a((\omega_i))=(01(011)^{z-1}010\cdots x_{6z+4})^\infty\preceq(01(011)^{z-1}01(011)^{z+1})^\infty$$ which is a contraction to (18). Hence we get $x_{3z+2}=1$ then $x_{3z+3}=0$ and $x_{3z+4}=1$.
							}
							\item[$\bullet$]
							{
								If $x_{3z+5}=0$, then $x_{3z+6}=1$ $$\sigma^a((\omega_i))=(01(011)^{z}0101\cdots x_{6z+4})^\infty.$$
								\begin{itemize}
									\item 
									{
										If $x_{3z+7}=0$, then $$\sigma^{a+3z+2}((\omega_i))=(01010x_{3z+8}\cdots x_{6z+4}01(011)^{z})^\infty\prec(01(011)^{z-1}01(011)^{z+1})^\infty$$which is a contraction to (18). Hence we have $x_{3z+7}=1$, then $x_{3z+8}=0$ and $x_{3z+9}=1$.
									}
									\item 
									{
										If $x_{3z+10}=0$, then $$\sigma^{a+3z+2}((\omega_i))=(01011010x_{3z+11}\cdots x_{6z+4}01(011)^{z})^\infty\prec(01(011)^{z-1}01(011)^{z+1})^\infty,$$which is a contraction. \\Hence we have $x_{3z+10}=1$, then $x_{3z+11}=0$ and $x_{3z+12}=1$.
									}\\
									Continue the same argument until $x_{6z+1}=0$, we will get a contraction. Hence we have $x_{6z+1}=1$ then $x_{6z+2}=0,x_{6z+3}=1$ and $$\sigma^a((\omega_i))=(01(011)^{z}01(011)^{z-1}01x_{6z+4})^\infty.$$
									\item
									{
										If $x_{6z+4}=1$, we have $\sigma^a((\omega_i))=(01(011)^{z})^\infty$ which is of the smallest period 3z+2 which is a contraction. Hence we have $x_{6z+4}=0$.\\
										However $\sigma^{a+6z+1}((\omega_i))=(01001(011)^{z}01(011)^{z-1})^\infty\prec(01(011)^{z-1}01(011)^{z+1})^\infty$ which is also a contraction to (18). 
									}
								\end{itemize}
								Hence we have $x_{3z+5}=1$ then $x_{3z+6}=0,x_{3z+7}=1$ and $$\sigma^a((\omega_i))=(01(011)^{z+1}01x_{3z+8}\cdots x_{6z+4})^\infty.$$ 
							}
							\item[$\bullet$]
							{
								If $x_{3z+8}=0$, then $$\sigma^a((\omega_i))=(01(011)^{z+1}010\cdots x_{6z+4})^\infty.$$We have
								$$\sigma^{a+3z+5}((\omega_i))=(010x_{3z+9}\cdots x_{6z+4}01(011)^{z+1})^\infty\preceq(01(011)^{z-1}01(011)^{z+1})^\infty$$ which is a contraction to (18). Thus we get $x_{3z+8}=1$ then $x_{3z+9}=0$ and $x_{3z+10}=1.$
							}
							\item[$\bullet$]
							{Continue the same argument until $x_{6z+2}=0$. We will get a contraction which implies $x_{6z+2}=1$ then $x_{6z+3}=0,x_{6z+4}=1$ and $$\sigma^a((\omega_i))=(01(011)^{2z}01)^\infty.$$
								However $\sigma^{a+6z+2}((\omega_i))=(0101(011)^{z})^\infty\prec(01(011)^{z-1}01(011)^{z+1})^\infty$, which is a contraction to (18).
							}
						\end{itemize} 
						Hence there is no such (6z+4)-period sequence $(\omega_i)\in \widetilde{K_\beta(t)}$ with $x_3=0$.
					}
					\item [\ding{193}]
					{
						If $x_3=1$, then $x_4=0,x_5=1$ and 
						$$\sigma^a((\omega_i))=(01101x_6\cdots x_{6z+4})^\infty.$$
						
						{
							If $x_6=0, \sigma^{a+3}((\omega_i))=:(y_i)=(010y_4\cdots y_{6z+1}011)^\infty$.\\
							Conduct the argument \ding{192} where $x_3=0$ for here $x_6=y_3=0$ and the following $x_9=0$ until $x_{6z}=0$, we will get a contraction which implies $x_{6z}=1$ then $x_{6z+1}=0$ and $$\sigma^a((\omega_i))=((011)^{2z}0011)^\infty.$$
							However $\sigma^{a+6z}((\omega_i))=(0011(011)^{2z})\prec(01(011)^{z-1}01(011)^{z+1})^\infty$ which is a contraction.
						}
						
						Hence there is no such (6z+4)-period sequence $(\omega_i)\in \widetilde{K_\beta(t)}$ with $x_3=1$. 
						
					}
					Hence there is no such (6z+4)-sequence $(\omega_i)\in \widetilde{K_\beta(t)}$ which is a contraction. \\Hence $\forall z\in \mathbb{N}_{\geq1},S_\beta(6z+4)\leq[(01(011)^{z-1}01(011)^{z+1})^\infty]_\beta$.  
				\end{itemize}
			}
			Combined 1) with 2), we obtain $\forall z\in \mathbb{N}_{\geq1},S_\beta(6z+4)=[(01(011)^{z-1}01(011)^{z+1})^\infty]_\beta$.
			
		}

	\end{proof}
\end{lem}

\begin{cor}
	By calculation, we obtain that \begin{equation*}
		S_\beta(3m+1)= \left\{ 
		\begin{aligned}
			&\frac{1+\beta+\beta^{3z+2}-\beta^{3z+1}-\beta^{6z}-\beta^{6z+2}}{(1-\beta^3)(\beta^{6z+1}-1)} , m=2z,z\in \mathbb{N}_{\geq1},\\
			&\frac{1+\beta+\beta^{3z+5}-\beta^{3z+4}-\beta^{6z+3}-\beta^{6z+5}}{(1-\beta^3)(\beta^{6z+4}-1)} , m=2z+1,z\in \mathbb{N}_{\geq1}.\\
		\end{aligned}  
		\right.        
	\end{equation*}
	The value of $S_\beta(p)$ is increasing to $\frac{\beta^{2}+1}{\beta^4-\beta}\approx 0.45626$ as p increasing to infinity.
	\begin{figure}[h]  
		\centering   
		\begin{minipage}{.5\textwidth}  
			\centering  
			\includegraphics[width=.8\linewidth]{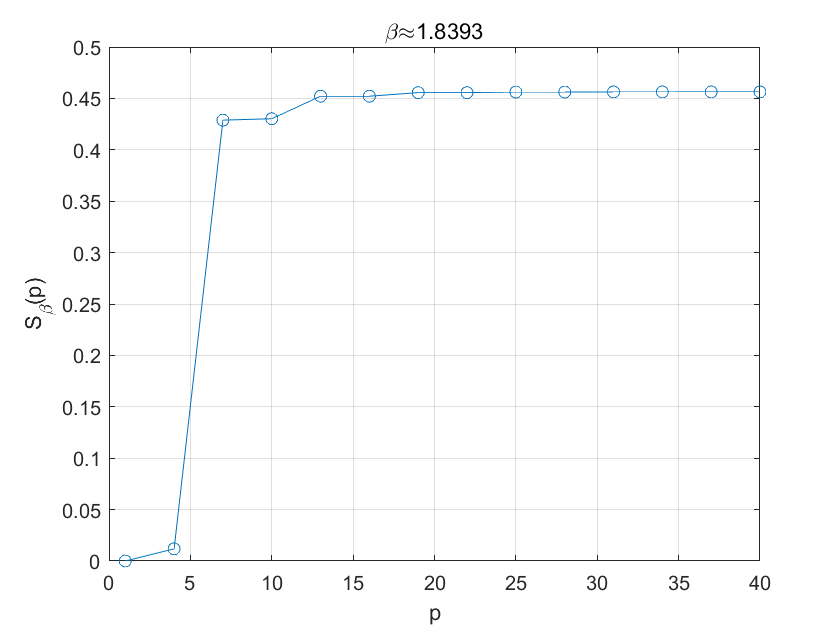}  
			\label{fig:image2}  
		\end{minipage}
		\tiny\caption{The value of $S_\beta(p)$ for $\beta\approx1.8393$ with the period $p=3m+1,\forall m\in\mathbb{N}$ from 1 to 40.}
	\end{figure} 
\end{cor}

\begin{lem}
	For $\beta=\phi_3,\ \delta(\beta)=(1^20)^\infty$, we have $$S_\beta(p)=[(01(011)^{z-1}01(011)^{z-1}01(011)^z)^\infty]_\beta,p=3m,m=3z,z \in \mathbb{N}_{\geq1}.$$
	\begin{proof}
		{
			(1)Firstly, we will prove that $$S_\beta(p)=S_\beta(9z)\geq [(01(011)^{z-1}01(011)^{z-1}01(011)^z)^\infty]_\beta=: t_*,\forall z\in \mathbb{N}_{\geq1}.$$
			We get that  $$\widetilde{K_\beta(t_*)}=\left\{(d_i)\in \left\{0,1\right\}^\mathbb{N} :(01(011)^{z-1}01(011)^{z-1}01(011)^z)^\infty\preceq\sigma^n((d_i))\prec(1^20)^\infty\right\}.$$
			Obviously the sequence $(01(011)^{z-1}01(011)^{z-1}01(011)^z)^\infty\in\widetilde{K_\beta(t_*)}$. \\
			That is $\widetilde{K_\beta(t_*)}$ contains a sequence of the smallest period 9z.\\
			By the definition of $S_\beta(9z)$, we have $S_\beta(9z)\geq t_*.$\\
			(2)Next to prove that $S_\beta(9z)\leq [(01(011)^{z-1}01(011)^{z-1}01(011)^z)^\infty]_\beta=: t_*$.\\
			Suppose on the contrary that $S_\beta(9z)>t_*$.\par
			Then there exists $t\in (t_*,S_\beta(9z))$, and then $\widetilde{K_\beta(t)}$ contains at least a periodic sequence $(\omega_i)$ of the smallest period 9z. And we have 
			\begin{equation}
				(01(011)^{z-1}01(011)^{z-1}01(011)^z)^\infty\prec(t_i)\preceq \sigma^j((\omega_i))\prec(1^20)^\infty.
			\end{equation}
			Since $(\omega_i)$ is of the smallest period 9z, there is at least one 1 among $\omega_1,\omega_2,\omega_3\cdots\omega_{9z}$.
			And since $\sigma^j((\omega_i))\prec(1^20)^\infty,\forall j \in \mathbb{N}$, there is neither consecutive three 1 nor consecutive 0 in $(\omega_i)$ which implies that there exists 0 among $\omega_1,\omega_2,\omega_3\cdots\omega_{9z}$.\par
			Therefore $\exists a \in \mathbb{N}$ such that $$\sigma^a((\omega_i))=(0x_2x_3\cdots x_{9z})^\infty.$$
			Since the inequality (19), we have $x_2=1$.
			\begin{itemize}
				\item [1)]
				{
					If $x_3=0.$\\
					Since $\sigma^j((\omega_i))\succ(01(011)^{z-1}01(011)^{z-1}01(011)^z)^\infty,\forall j\in \mathbb{N}$, we have 
					$$\sigma^a((\omega_i))=(01(011)^{z-1}01x_{3z+2}\cdots x_{9z})^\infty.$$
					\begin{itemize}
						\item [(I)]
						{
							If $x_{3z+2}=0.$\\
							Since $\sigma^j((\omega_i))\succ(01(011)^{z-1}01(011)^{z-1}01(011)^z)^\infty,\forall j\in \mathbb{N}$, we have 
							$$\sigma^a((\omega_i))=(01(011)^{z-1}01(011)^{z-1}01x_{6z+1}\cdots x_{9z})^\infty.$$
							\begin{itemize}
								\item [(i)]
								{
									If $x_{6z+1}=0$ then $$\sigma^a((\omega_i))=(01(011)^{z-1}01(011)^{z-1}010\cdots x_{9z})^\infty\preceq(01(011)^{z-1}01(011)^{z-1}01(011)^z)^\infty$$ which is a contraction to (19). \\
									Hence we have $x_{6z+1}=1$ then $x_{6z+2}=0,x_{6z+3}=1$ and $$\sigma^a((\omega_i))=(01(011)^{z-1}01(011)^{z-1}01101x_{6z+4}\cdots x_{9z})^\infty.$$
								}
								\item [(ii)]
								{
									If $x_{6z+4}=0$ then $\sigma^a((\omega_i))=(01(011)^{z-1}01(011)^{z-1}011010\cdots x_{9z})^\infty$ and
									
									\begin{align}
										\sigma^{a+6z+1}((\omega_i))&=(010x_{6z+5}\cdots x_{9z}01(011)^{z-1}01(011)^{z})\nonumber\\
										&\preceq(01(011)^{z-1}01(011)^{z-1}01(011)^z)^\infty \nonumber
									\end{align}
									which is a contraction to (19). \\
									Hence we have $x_{6z+4}=1$ then $x_{6z+5}=0,x_{6z+6}=1$ and $$\sigma^a((\omega_i))=(01(011)^{z-1}01(011)^{z+1}01x_{6z+7}\cdots x_{9z})^\infty.$$
								}
								\item [(iii)]
								{
									If $x_{6z+7}=0$ then $\sigma^a((\omega_i))=(01(011)^{z-1}01(011)^{z+1}010\cdots x_{9z})^\infty$ and
									\begin{align}
										\sigma^{a+6z+4}((\omega_i))&=(010x_{6z+8}\cdots x_{9z}01(011)^{z-1}01(011)^{z+1})^\infty\nonumber\\
										&\preceq(01(011)^{z-2}01(011)^{z-1}01(011)^{z+1})^\infty \nonumber\\
										&\prec(01(011)^{z-1}01(011)^{z-1}01(011)^z)^\infty \nonumber
									\end{align}
									which is a contraction to (19).\\
									Hence we have $x_{6z+7}=1$ then $x_{6z+8}=0,x_{6z+9}=1$ and $$\sigma^a((\omega_i))=(01(011)^{z-1}01(011)^{z+2}01x_{6z+10}\cdots x_{9z})^\infty.$$
								}
								Continue the same argument until $x_{9z-2}=0$, we will get a contraction which implies $x_{9z-2}=1,x_{9z-1}=0,x_{9z}=1$ and $$\sigma^a((\omega_i))=(01(011)^{z-1}01(011)^{2z-1}01)^\infty$$
								Hence we have 
								$$\sigma^{a+9z-2}((\omega_i))=(0101(011)^{z-1}01(011)^{2z-1})^\infty\prec(01(011)^{z-1}01(011)^{z-1}01(011)^z)^\infty$$ which is a contraction to (19).
							\end{itemize}
							Hence we have $x_{3z+2}=1$ then $x_{3z+3}=0,x_{3z+4}=1$ and $$\sigma^a((\omega_i))=(01(011)^z01x_{3z+5}\cdots x_{9z})^\infty.$$
							
						}
						\item [(II)]
						{
							If $x_{3z+5}=0$, then we have $x_{3z+6}=1$ and $$\sigma^{a+3z+2}((\omega_i))=(0101x_{3z+7}\cdots x_{9z}01(011)^z)^\infty.$$
							Since $\sigma^{a+3z+2}((\omega_i))\succ(01(011)^{z-1}01(011)^{z-1}01(011)^z)^\infty$, we have $$\sigma^{a+3z+2}((\omega_i))=(01(011)^{z-1}01x_{6z+4}\cdots x_{9z}01(011)^z)^\infty.$$
							\begin{itemize}
								\item [(i)]
								{
									If $x_{6z+4}=0$, then $x_{6z+5}=1$ 
									\begin{align}
										\sigma^{a+3z+2}((\omega_i))&=(01(011)^{z-1}0101x_{6z+6}\cdots x_{9z}01(011)^z)^\infty\nonumber\\
										&\preceq(01(011)^{z-1}01(011)^{z-1}01(011)^z)^\infty \nonumber
									\end{align}
									which is a contraction to (19). \\
									Hence we have $x_{6z+4}=1$ then $x_{6z+5}=0,x_{6z+6}=1$ and $$\sigma^{a+3z+2}((\omega_i))=(01(011)^{z}01x_{6z+7}\cdots x_{9z}01(011)^z)^\infty.$$
								}
								\item{(ii)}
								{
									If $x_{6z+7}=0$, then we obtain a contraction that 
									\begin{align}
										\sigma^{a+3z+2+3z+2}((\omega_i))&=\sigma^{a+6z+4}((\omega_i)) \nonumber\\
										&=(010x_{6z+8}\cdots x_{9z}01(011)^z01(011)^z)^\infty \nonumber\\
										&\preceq(01(011)^{z-2}01(011)^{z}01(011)^z)^\infty \nonumber\\
										&\prec(01(011)^{z-1}01(011)^{z-1}01(011)^z)^\infty.\nonumber
									\end{align}
									Hence we have $x_{6z+7}=1$ then $x_{6z+8}=0$ and $x_{6z+9}=1.$\\	 
								}
								Continue the same argument until $x_{9z-2}=0$, we will get a contraction which implies that $x_{9z-2}=1$ then $x_{9z-1}=0,x_{9z}=1$ and $$\sigma^{a+3z+2}((\omega_i))=(01(011)^{2z-2}0101(011)^z)^\infty.$$
								However we find that
								\begin{align}
									\sigma^{a+3z+2+6z-4}((\omega_i))&=\sigma^{a+9z-2}((\omega_i)) \nonumber\\
									&=(0101(011)^{z}01(011)^{2z-2})^\infty \nonumber\\
									&\prec(01(011)^{z-1}01(011)^{z-1}01(011)^z)^\infty .\nonumber
								\end{align}
								which is a contraction to (19).
							\end{itemize}
							Hence we have $x_{3z+5}=1$ then $x_{3z+6}=0,x_{3z+7}=1$ and $$\sigma^a((\omega_i))=(01(011)^{z+1}01x_{3z+8}\cdots x_{9z})^\infty.$$
						}
						\item [(III)]
						{
							Continue the same argument as $x_{3z+2}=0$ and $x_{3z+5}=0$ for $x_{3z+8}=0$, we will get a contraction which implies that $x_{3z+8}=1$ and then $x_{3z+9}=0$ and $x_{3z+10}=1$.
						}
						Continue the argument until $x_{9z-4}=0$, we will get a contraction which implies that $x_{9z-4}=1$ then $x_{9z-3}=0,x_{9z-2}=1$ and $$\sigma^a((\omega_i))=(01(011)^{3z-2}01x_{9z-1}x_{9z})^\infty.$$
						\item [(IV)]
						{
							If $x_{9z-1}=0$, then we have $x_{9z}=1$. But $$\sigma^{a+9z-4}((\omega_i))=(010101(011)^{3z-2})^\infty\prec(01(011)^{z-1}01(011)^{z-1}01(011)^z)^\infty$$ which is a contraction to (19). Hence we have $x_{9z-1}=1$ and then $x_{9z}=0$.
							However we have $$\sigma^{a+9z-1}((\omega_i))=(001(011)^{3z-1})^\infty\prec(01(011)^{z-1}01(011)^{z-1}01(011)^z)^\infty$$ which is also a contraction to (19).
						}
					\end{itemize}
					Hence we have $x_3=1$ then $x_4=0,x_5=1$ and $$\sigma^a((\omega_i))=(01101x_6\cdots x_{9z})^\infty.$$
					And $\sigma^{a+3}((\omega_i))=(01x_6\cdots x_{9z}011)^\infty=:(y_i)=(01y_3\cdots y_{9z-3}011)^\infty.$
				}
				\item [2)]
				{
					If $x_6=y_3=0.$\\
					Since $\sigma^{a+3}((\omega_i))=(01y_3\cdots y_{9z-3}011)^\infty\succ(01(011)^{z-1}01(011)^{z-1}01(011)^z)^\infty,$ we have 
					$$(y_i)=(01(011)^{z-1}01y_{3z+2}\cdots y_{9z-3}011)^\infty$$ which is the same situation as the step \ding{192} where $\sigma^a((\omega_i))=(01(011)^{z-1}01x_{3z+2}\cdots x_{9z})^\infty$.
					
				}
				\item [2)]
				{
					Hence we conduct the same argument as \ding{192} and \ding{193} until $x_{9z-3}=y_{9z-6}=0$ which leads to a contraction which implies that $y_{9z-6}=1$ then $y_{9z-5}=0,y_{9z-4}=1$ and $$(y_i)=((011)^{3z-2}01y_{9z-3}011)^\infty.$$
					If $y_{9z-3}=1$, then $(y_i)=((011)^{3z})^\infty$ which is of the smallest period 3 which is a contraction. Hence we have $y_{9z-3}=0$.
					However we find that $$\sigma^{a+9z-6}(((y_i))=\sigma^{a+9z-3}((\omega_i))=(010011(011)^{3z-2})^\infty\prec(01(011)^{z-1}01(011)^{z-1}01(011)^z)^\infty.$$ which is also a contraction to (19).
				}
			\end{itemize}
			Hence there is no such 9z-period sequence $(\omega_i)\in \widetilde{K_\beta(t)}$ which is a contraction. Hence we have $$S_\beta(9z)\leq [(01(011)^{z-1}01(011)^{z-1}01(011)^z)^\infty]_\beta.$$
			Combined 1) and 2) we get that $$S_\beta(9z)= [(01(011)^{z-1}01(011)^{z-1}01(011)^z)^\infty]_\beta,\forall z \in \mathbb{N}_{\geq1}.$$
		}
	\end{proof}
\end{lem}

\begin{lem}
	For $\beta=\phi_3,\ \delta(\beta)=(1^20)^\infty$, we have $$S_\beta(p)=[(01(011)^{z-1}01(011)^z01(011)^z)^\infty]_\beta ,p=3m,m=3z+1,z\in \mathbb{N}_{\geq1}.$$
	\begin{proof}
		(1)Firstly, we will prove that $$S_\beta(p)=S_\beta(9z+3)\geq [(01(011)^{z-1}01(011)^z01(011)^z)^\infty]_\beta=: t_*,\forall z\in \mathbb{N}_{\geq1}.$$
		We get that  $$\widetilde{K_\beta(t_*)}=\left\{(d_i)\in \left\{0,1\right\}^\mathbb{N} :(01(011)^{z-1}01(011)^z01(011)^z)^\infty\preceq\sigma^n((d_i))\prec(1^20)^\infty\right\}.$$
		Obviously the sequence $$(01(011)^{z-1}01(011)^z01(011)^z)^\infty\in\widetilde{K_\beta(t_*)}.$$ That is $\widetilde{K_\beta(t_*)}$ contains a sequence of the smallest period 9z+3.\\
		By the definition of $S_\beta(9z+3)$, we have $S_\beta(9z+3)\geq t_*.$\\
		(2)Next to prove that $S_\beta(9z+3)\leq [(01(011)^{z-1}01(011)^z01(011)^z)^\infty]_\beta=: t_*$.\\
		Suppose on the contrary that $S_\beta(9z+3)>t_*$.
		Then there exists $t\in (t_*,S_\beta(9z+3))$, and then $\widetilde{K_\beta(t)}$ contains at least a periodic sequence $(\omega_i)$ of the smallest period 9z+3. And we have 
		\begin{equation}
			(01(011)^{z-1}01(011)^z01(011)^z)^\infty\prec(t_i)\preceq \sigma^j((\omega_i))\prec(1^20)^\infty.
		\end{equation}
		Since $(\omega_i)$ is of the smallest period 9z+3, there is at least one 1 among $\omega_1,\omega_2,\omega_3\cdots\omega_{9z+3}$.\par
		And since $\sigma^j((\omega_i))\prec(1^20)^\infty,\forall j \in \mathbb{N}$, there is neither consecutive three 1 nor consecutive 0 in $(\omega_i)$  which implies that there exists 0 among $\omega_1,\omega_2,\omega_3\cdots\omega_{9z+3}$.\par
		Therefore $\exists a \in \mathbb{N}$ such that $$\sigma^a((\omega_i))=(0x_2x_3\cdots x_{9z+3})^\infty.$$
		Since the inequality (20), we have $x_2=1$.
		\begin{itemize}
			\item [1)]
			{
				If $x_3=0.$\\
				Since $\sigma^j((\omega_i))\succ(01(011)^{z-1}01(011)^z01(011)^z)^\infty,\forall j\in \mathbb{N}$, we have 
				$$\sigma^a((\omega_i))=(01(011)^{z-1}01x_{3z+2}\cdots x_{9z+3})^\infty.$$
				\begin{itemize}
						\item [(I)]
					{
						If $x_{3z+2}=0.$\\
						Since $\sigma^j((\omega_i))\succ(01(011)^{z-1}01(011)^z01(011)^z)^\infty,\forall j\in \mathbb{N}$, we have 
						$$\sigma^a((\omega_i))=(01(011)^{z-1}01(011)^z01x_{6z+4}\cdots x_{9z+3})^\infty.$$
						\begin{itemize}
							\item [(i)]
							{
								If $x_{6z+4}=0$ then $$\sigma^a((\omega_i))=(01(011)^{z-1}01(011)^z010\cdots x_{9z+3})^\infty\preceq(01(011)^{z-1}01(011)^z01(011)^z)^\infty$$ which is a contraction to (20).\\
								Hence we have $x_{6z+4}=1$ then $x_{6z+5}=0,x_{6z+6}=1$ and $$\sigma^a((\omega_i))=(01(011)^{z-1}01(011)^{z+1}01x_{6z+7}\cdots x_{9z+3})^\infty.$$
							}
							\item [(ii)]
							{
								If $x_{6z+7}=0$ then $\sigma^a((\omega_i))=(01(011)^{z-1}01(011)^{z+1}010\cdots x_{9z+3})^\infty$ and
								\begin{align}
									\sigma^{a+6z+4}((\omega_i))&=(010x_{6z+8}\cdots x_{9z+3}01(011)^{z-1}01(011)^{z+1})^\infty\nonumber\\
									&\preceq(01(011)^{z-1}01(011)^{z-1}01(011)^{z+1})^\infty \nonumber\\
									&\prec(01(011)^{z-1}01(011)^z01(011)^z)^\infty \nonumber
								\end{align}
								which is a contraction.\\
								Hence we have $x_{6z+7}=1$ then $x_{6z+8}=0,x_{6z+9}=1$ and $$\sigma^a((\omega_i))=(01(011)^{z-1}01(011)^{z+2}01x_{6z+10}\cdots x_{9z+3})^\infty.$$
							}
							\item [(iii)]
							{
								If $x_{6z+10}=0$ then $\sigma^a((\omega_i))=(01(011)^{z-1}01(011)^{z+2}010\cdots x_{9z+3})^\infty$ and
								\begin{align}
									\sigma^{a+6z+7}((\omega_i))&=(010x_{6z+11}\cdots x_{9z+3}01(011)^{z-1}01(011)^{z+2})^\infty\nonumber\\
									&\preceq(01(011)^{z-2}01(011)^{z-1}01(011)^{z+2})^\infty \nonumber\\
									&\prec(01(011)^{z-1}01(011)^z01(011)^z)^\infty \nonumber
								\end{align}
								which is a contraction. \\
								Hence we have $x_{6z+10}=1$ then $x_{6z+11}=0,x_{6z+12}=1$ and $$\sigma^a((\omega_i))=(01(011)^{z-1}01(011)^{z+3}01x_{6z+13}\cdots x_{9z+3})^\infty.$$
							}
							Continue the same argument until $x_{9z+1}=0$, we will get a contraction which implies $x_{9z+1}=1,x_{9z+2}=0,x_{9z+3}=1$ and $$\sigma^a((\omega_i))=(01(011)^{z-1}01(011)^{2z}01)^\infty.$$
							Hence we have 
							$$\sigma^{a+9z+1}((\omega_i))=(0101(011)^{z-1}01(011)^{2z})^\infty\prec(01(011)^{z-1}01(011)^z01(011)^z)^\infty$$ which is a contraction.
						\end{itemize}
						Hence we have $x_{3z+2}=1$ then $x_{3z+3}=0,x_{3z+4}=1$ and $$\sigma^a((\omega_i))=(01(011)^z01x_{3z+5}\cdots x_{9z+3})^\infty.$$
						
					}
					\item [(II)]
					{
						If $x_{3z+5}=0$, then we have $x_{3z+6}=1$ and $$\sigma^{a+3z+2}((\omega_i))=(0101x_{3z+7}\cdots x_{9z+3}01(011)^z)^\infty.$$
						Since $\sigma^{a+3z+2}((\omega_i))\succ(01(011)^{z-1}01(011)^z01(011)^z)^\infty$, we have $$\sigma^{a+3z+2}((\omega_i))=(01(011)^{z-1}01x_{6z+4}\cdots x_{9z+3}01(011)^z)^\infty.$$
						\begin{itemize}
							\item [(i)]
							{
								If $x_{6z+4}=0$, then $x_{6z+5}=1$ $$\sigma^{a+3z+2}((\omega_i))=(01(011)^{z-1}0101\cdots x_{9z+3}01(011)^z)^\infty\preceq(01(011)^{z-1}01(011)^{z}01(011)^z)^\infty$$ which is a contraction to (20). \\
								Hence we have $x_{6z+4}=1$ then $x_{6z+5}=0,x_{6z+6}=1$ and $$\sigma^{a+3z+2}((\omega_i))=(01(011)^{z}01x_{6z+7}\cdots x_{9z+3}01(011)^z)^\infty.$$
							}
							\item[(ii)]
							{
								If $x_{6z+7}=0$, then we obtain a contraction that 
								\begin{align}
									\sigma^{a+3z+2+3z+2}((\omega_i))&=\sigma^{a+6z+4}((\omega_i)) \nonumber\\
									&=(010x_{6z+8}\cdots x_{9z+3}01(011)^z01(011)^z)^\infty \nonumber\\
									&\preceq(01(011)^{z-1}01(011)^{z}01(011)^z)^\infty \nonumber
								\end{align}
								Hence we have $x_{6z+7}=1$ then $x_{6z+8}=0$ and $x_{6z+9}=1$.\\	 
							}
							Continue the same argument until $x_{9z+1}=0$, we will get a contraction which implies that $x_{9z+1}=1$ then $x_{9z+2}=0,x_{9z+3}=1$ and $$\sigma^{a+3z+2}((\omega_i))=(01(011)^{2z-1}0101(011)^z)^\infty.$$
							However we find that
							\begin{align}
								\sigma^{a+3z+2+6z-1}((\omega_i))&=\sigma^{a+9z+1}((\omega_i)) \nonumber\\
								&=(0101(011)^{z}01(011)^{2z-1})^\infty \nonumber\\
								&\prec(01(011)^{z-1}01(011)^z01(011)^z)^\infty \nonumber
							\end{align}
							which is a contraction to (20).
						\end{itemize}
						Hence we have $x_{3z+5}=1$ then $x_{3z+6}=0,x_{3z+7}=1$ and $$\sigma^a((\omega_i))=(01(011)^{z+1}01x_{3z+8}\cdots x_{9z+3})^\infty.$$
					}
					\item [(III)]
					{
						Continue the same argument as $x_{3z+2}=0$ and $x_{3z+5}=0$ for $x_{3z+8}=0$, we will get a contraction which implies that $x_{3z+8}=1$ and then $x_{3z+9}=0$ and $x_{3z+10}=1$.
					}
					Continue the argument until $x_{9z-1}=0$, we will get a contraction which implies that $x_{9z-1}=1$ then $x_{9z}=0,x_{9z+1}=1$ and $$\sigma^a((\omega_i))=(01(011)^{3z-1}01x_{9z+2}x_{9z+3})^\infty.$$
					\item [(IV)]
					{
						If $x_{9z+2}=0$, then we have $x_{9z+3}=1$.\\
						But $\sigma^{a+9z-1}((\omega_i))=(010101(011)^{3z-1})^\infty<(01(011)^{z-1}01(011)^z01(011)^z)^\infty$ which is a contraction. Hence we have $x_{9z+2}=1$ and then $x_{9z+3}=0$.\\
						However we have $\sigma^{a+9z+2}((\omega_i))=(001(011)^{3z})^\infty\prec(01(011)^{z-1}01(011)^z01(011)^z)^\infty$ which is also a contraction.
					}
				\end{itemize}
				Hence we have $x_3=1$ then $x_4=0,x_5=1$ and $$\sigma^a((\omega_i))=(01101x_6\cdots x_{9z+3})^\infty.$$
				And $\sigma^{a+3}((\omega_i))=(01x_6\cdots x_{9z+3}011)^\infty=:{(y_i)}=(01y_3\cdots y_{9z}011)^\infty.$
			}
			\item [2)]
			{
				If $x_6=y_3=0.$\\
				Since $\sigma^{a+3}((\omega_i))=(01y_3\cdots y_{9z}011)^\infty\succ(01(011)^{z-1}01(011)^z01(011)^z)^\infty,$ we have 
				$$(y_i)=(01(011)^{z-1}01y_{3z+2}\cdots y_{9z}011)^\infty$$ which is the same situation as the step \ding{192} where $\sigma^a((\omega_i))=(01(011)^{z-1}01x_{3z+2}\cdots x_{9z+3})^\infty$.\\
				Therefore we can get a contraction which leads to $$x_6=y_3=1,x_7=y_4=0,x_8=y_5=1.$$
				
			}
			\item [3)]
			{
				Hence we conduct the same argument as \ding{192} and \ding{193} until $x_{9z}=y_{9z-3}=0$ which leads to a contraction which implies that $y_{9z-3}=1$ then $y_{9z-2}=0,y_{9z-1}=1$ and $$(y_i)=((011)^{3z-1}01y_{9z}011)^\infty.$$
				If $y_{9z}=1$, then $(y_i)=((011)^{3z+1})^\infty$ which is of the smallest period 3 which is a contraction. Hence we have $y_{9z}=0$.However we find that $$\sigma^{a+9z-3}(((y_i))=\sigma^{a+9z}((\omega_i))=(010011(011)^{3z-1})^\infty\prec(01(011)^{z-1}01(011)^z01(011)^z)^\infty,$$ which is also a contraction to (20).
			}
		\end{itemize}
		Hence there is no such (9z+3)-period sequence $(\omega_i)\in \widetilde{K_\beta(t)}$ which is a contraction. Hence we have $$S_\beta(9z+3)\leq [(01(011)^{z-1}01(011)^z01(011)^z)^\infty]_\beta.$$
		Combined 1) and 2) we get that $$S_\beta(9z+3)= [(01(011)^{z-1}01(011)^z01(011)^z)^\infty]_\beta,\forall z \in \mathbb{N}_{\geq1}.$$
		
	\end{proof}
\end{lem}

\begin{them}
	{
		For $\beta=\phi_3,\ \delta(\beta)=(1^20)^\infty$, we have $S_\beta(1)=0,S_\beta(3)=[(001)^\infty]_\beta,S_\beta(4)=[(0011)^\infty]_\beta,S_\beta(6)=[(001101)^\infty]_\beta$ and
		\begin{equation*}
			S_\beta(p)= \left\{ 
			\begin{aligned}
				&[01(011)^{m})^\infty]_\beta,p=3m+2,\forall m \in \mathbb{N},\\	&[(01(011)^{z-1}01(011)^z)^\infty]_\beta ,p=3m+1,m=2z,z\in \mathbb{N}_{\geq1},\\
				&[(01(011)^{z-1}01(011)^{z+1})^\infty]_\beta ,p=3m+1,m=2z+1,z\in \mathbb{N}_{\geq1}.\\
				&[(01(011)^{z-1}01(011)^{z-1}01(011)^z)^\infty]_\beta,p=3m,m=3z,z \in \mathbb{N}_{\geq1},\\	&[(01(011)^{z-1}01(011)^z01(011)^z)^\infty]_\beta ,p=3m,m=3z+1,z\in \mathbb{N}_{\geq1},\\
				&[(01(011)^{z-1}01(011)^{z+1}01(011)^z)^\infty]_\beta,p=3m,m=3z+2,z\in \mathbb{N}_{\geq1}.\\
			\end{aligned}  
			\right.        
		\end{equation*}
	}
	\begin{proof}
		Since we have got many results for the different choice of the value of $p$. Now we need to prove the left one result above
		$$S_\beta(p)=[(01(011)^{z-1}01(011)^{z+1}01(011)^z)]_\beta^\infty,p=3m,m=3z+2,z \in \mathbb{N}_{\geq1}.$$
		{
			(1)Firstly, we will prove that $$S_\beta(p)=S_\beta(9z+6)\geq [(01(011)^{z-1}01(011)^{z+1}01(011)^z)^\infty]_\beta=: t_*,\forall z\in \mathbb{N}_{\geq1}.$$
			We get that  $$\widetilde{K_\beta(t_*)}=\left\{(d_i)\in \left\{0,1\right\}^\mathbb{N} :(01(011)^{z-1}01(011)^{z+1}01(011)^z)^\infty\preceq\sigma^n((d_i))\prec(1^20)^\infty\right\}.$$
			Obviously the sequence $(01(011)^{z-1}01(011)^{z+1}01(011)^z)^\infty\in\widetilde{K_\beta(t_*)}$. That is $\widetilde{K_\beta(t_*)}$ contains a sequence of the smallest period 9z+6.\\
			By the definition of $S_\beta(9z+6)$, we have $S_\beta(9z+6)\geq t_*.$\\
			(2)Next to prove that $S_\beta(9z+6)\leq [(01(011)^{z-1}01(011)^{z+1}01(011)^z)^\infty]_\beta=: t_*$.\\
			Suppose on the contrary that $S_\beta(9z+6)>t_*$.\par
			Then there exists $t\in (t_*,S_\beta(9z+6))$, and then $\widetilde{K_\beta(t)}$ contains at least a periodic sequence $(\omega_i)$ of the smallest period 9z+6. And we have 
			\begin{equation}
				(01(011)^{z-1}01(011)^{z+1}01(011)^z)^\infty\prec(t_i)\preceq \sigma^j((\omega_i))\prec(1^20)^\infty.
			\end{equation}
			Since $(\omega_i)$ is of the smallest period 9z+6, there is at least one 1 among $\omega_1,\omega_2,\omega_3\cdots\omega_{9z+6}$.\par
			And since $\sigma^j((\omega_i))\prec(1^20)^\infty,\forall j \in \mathbb{N}$, there is neither consecutive three 1 nor consecutive 0 in $(\omega_i)$  which implies that there exists 0 among $\omega_1,\omega_2,\omega_3\cdots\omega_{9z+6}$.
			Therefore $\exists a \in \mathbb{N}$ such that $$\sigma^a((\omega_i))=(0x_2x_3\cdots x_{9z+6})^\infty.$$
			Since the inequality (21), we have $x_2=1$.
			\begin{itemize}
				\item [1)]
				{
					If $x_3=0.$\\
					Since $\sigma^j((\omega_i))\succ(01(011)^{z-1}01(011)^{z+1}01(011)^z)^\infty,\forall j\in \mathbb{N}$, we have 
					$$\sigma^a((\omega_i))=(01(011)^{z-1}01x_{3z+2}\cdots x_{9z+6})^\infty.$$
					\begin{itemize}
						\item [\ding{192}]
						{
							If $x_{3z+2}=0.$\\
							Since $\sigma^j((\omega_i))\succ(01(011)^{z-1}01(011)^{z+1}01(011)^z)^\infty,\forall j\in \mathbb{N}$, we have 
							$$\sigma^a((\omega_i))=(01(011)^{z-1}01(011)^{z+1}01x_{6z+7}\cdots x_{9z+6})^\infty.$$
							\begin{itemize}
								\item [(I)]
								{
									If $x_{6z+7}=0$ then $$\sigma^a((\omega_i))=(01(011)^{z-1}01(011)^{z+1}010\cdots x_{9z+6})^\infty\preceq(01(011)^{z-1}01(011)^{z+1}01(011)^z)^\infty$$ which is a contraction to (21). \\Hence we have $x_{6z+7}=1$ then $x_{6z+8}=0,x_{6z+9}=1$ and $$\sigma^a((\omega_i))=(01(011)^{z-1}01(011)^{z+2}01x_{6z+10}\cdots x_{9z+6})^\infty.$$
								}
								\item [(II)]
								{
									If $x_{6z+10}=0$ then $\sigma^a((\omega_i))=(01(011)^{z-1}01(011)^{z+2}010\cdots x_{9z+6})^\infty$ and
									\begin{align}
										\sigma^{a+6z+7}((\omega_i))&=(010x_{6z+11}\cdots x_{9z+6}01(011)^{z-1}01(011)^{z+2})^\infty\nonumber\\
										&\preceq(01(011)^{z-1}01(011)^{z-1}01(011)^{z+2})^\infty \nonumber\\
										&\prec(01(011)^{z-1}01(011)^{z+1}01(011)^z)^\infty \nonumber
									\end{align}
									which is a contraction. \\Hence we have $x_{6z+10}=1$ then $x_{6z+11}=0,x_{6z+12}=1$ and $$\sigma^a((\omega_i))=(01(011)^{z-1}01(011)^{z+3}01x_{6z+13}\cdots x_{9z+6})^\infty.$$
								}
								\item [(III)]
								{
									If $x_{6z+13}=0$ then $\sigma^a((\omega_i))=(01(011)^{z-1}01(011)^{z+3}010\cdots x_{9z+6})^\infty$ and
									\begin{align}
										\sigma^{a+6z+10}((\omega_i))&=(010x_{6z+14}\cdots x_{9z+6}01(011)^{z-1}01(011)^{z+3})^\infty\nonumber\\
										&\preceq(01(011)^{z-2}01(011)^{z-1}01(011)^{z+3})^\infty \nonumber\\
										&\prec(01(011)^{z-1}01(011)^{z+1}01(011)^z)^\infty\nonumber
									\end{align}
									which is a contraction. Hence we have $x_{6z+13}=1$ then $x_{6z+14}=0,x_{6z+15}=1$ and $\sigma^a((\omega_i))=(01(011)^{z-1}01(011)^{z+4}01x_{6z+16}\cdots x_{9z+6})^\infty.$\\
								}
								Continue the same argument until $x_{9z+4}=0$, we will get a contraction which implies $x_{9z+4}=1,x_{9z+5}=0,x_{9z+6}=1$ and $$\sigma^a((\omega_i))=(01(011)^{z-1}01(011)^{2z+1}01)^\infty$$
								Hence we have 
								$$\sigma^{a+9z+4}((\omega_i))=(0101(011)^{z-1}01(011)^{2z+1})^\infty\prec(01(011)^{z-1}01(011)^{z+1}01(011)^z)^\infty$$ which is a contraction.
							\end{itemize}
							Hence we have $x_{3z+2}=1$ then $x_{3z+3}=0,x_{3z+4}=1$ and $$\sigma^a((\omega_i))=(01(011)^z01x_{3z+5}\cdots x_{9z+6})^\infty.$$
							
						}
						\item [\ding{193}]
						{
							If $x_{3z+5}=0$, then we have $x_{3z+6}=1$ and $$\sigma^{a+3z+2}((\omega_i))=(0101x_{3z+7}\cdots x_{9z+6}01(011)^z)^\infty.$$
							Since $\sigma^{a+3z+2}((\omega_i))\succ(01(011)^{z-1}01(011)^{z+1}01(011)^z)^\infty$, we have $$\sigma^{a+3z+2}((\omega_i))=(01(011)^{z-1}01x_{6z+4}\cdots x_{9z+6}01(011)^z)^\infty$$
							\begin{itemize}
								\item [(I)]
								{
									If $x_{6z+4}=0$, then $x_{6z+5}=1$
									\begin{align}
										\sigma^{a+3z+2}((\omega_i))&=(01(011)^{z-1}0101\cdots x_{9z+6}01(011)^z)^\infty\nonumber\\
										&\preceq(01(011)^{z-1}01(011)^{z+1}01(011)^z)^\infty \nonumber
									\end{align}
									  which is a contraction to (21). \\Hence we have $x_{6z+4}=1$ then $x_{6z+5}=0,x_{6z+6}=1$ and $$\sigma^{a+3z+2}((\omega_i))=(01(011)^{z}01x_{6z+7}\cdots x_{9z+6}01(011)^z)^\infty.$$
								}
								\item[(II)]
								{
									If $x_{6z+7}=0$, then $x_{6z+8}=1$.
									\begin{itemize}
										\item [(i)]
										{
											If $x_{6z+9}=0$, then 
											\begin{align}
												\sigma^{a+6z+4}((\omega_i))&=(01010x_{6z+10}\cdots x_{9z+6}01(011)^z)\nonumber\\
												&\prec(01(011)^{z-1}01(011)^{z+1}01(011)^z)^\infty \nonumber
											\end{align}
											  which is a contraction. \\Hence we have $x_{6z+9}=1$ then $x_{6z+10}=0$ and $x_{6z+11}=1$.
											
										}
										\item [(ii)]
										{
											If $x_{6z+12}=0$, then 
											\begin{align}
												\sigma^{a+6z+4}((\omega_i))&=(01011010x_{6z+13}\cdots x_{9z+6}01(011)^z)\nonumber\\
												&\prec(01(011)^{z-1}01(011)^{z+1}01(011)^z)^\infty\nonumber
											\end{align}
											which is a contraction.\\ Hence we have $x_{6z+12}=1$ then $x_{6z+13}=0$ and $x_{6z+14}=1$.
										}
										\item [(iii)]
										{
											Continue the same argument until $x_{9z+3}=0$, we will get a contraction which implies $x_{9z+3}=1$ then $x_{9z+4}=0,x_{9z+5}=1$ and $$\sigma^{a+3z+2}((\omega_i))=(01(011)^{z}01(011)^{z-1}01x_{9z+6}01(011)^{z})^\infty.$$
											If $x_{9z+6}=1,\sigma^{a+3z+2}((\omega_i))=(01(011)^z)^\infty$ which is of the smallest period 3z+2 which is a contraction. Hence $x_{9z+6}=0$. However we find that
											\begin{align}
												\sigma^{a+9z+5}((\omega_i))&=(001(011)^z01(011)^{z}01(011)^{z-1}01)^\infty\nonumber\\
												&\prec(01(011)^{z-1}01(011)^{z+1}01(011)^z)^\infty\nonumber
											\end{align}
											  which is also a contraction.		  	
										}
									\end{itemize}
									Hence we have $x_{6z+7}=1$, then $x_{6z+8}=0,x_{6z+9}=1$ and $$\sigma^a((\omega_i))=(01(011)^{z}01(011)^{z+1}01x_{6z+10}\cdots x_{9z+6})^\infty$$
									\item[(III)]
									{
										If $x_{6z+10}=0$, then $x_{6z+11}=1$
										\begin{align}
											\sigma^{a+6z+7}((\omega_i))&=(010x_{6z+11}\cdots x_{9z+6}(011)^{z}01(011)^{z+1})^\infty\nonumber\\
											&\preceq(01(011)^{z-1}01(011)^z01(011)^{z+1})^\infty \nonumber\\
											&\prec(01(011)^{z-1}01(011)^{z+1}01(011)^z)^\infty \nonumber
										\end{align} 
										which is a contraction. Hence we have $x_{6z+10}=1$ then $x_{6z+11}=0,x_{6z+12}=1$ and $$\sigma^{a+3z+2}((\omega_i))=(01(011)^{z+2}01x_{6z+13}\cdots x_{9z+6}01(011)^z)^\infty.$$
									}
									\item [(IV)]
									{
										If $x_{6z+13}=0$, then we obtain a contraction that 
										\begin{align}
											\sigma^{a+3z+2+3z+8}((\omega_i))&=\sigma^{a+6z+10}((\omega_i)) \nonumber\\
											&=(010x_{6z+13}\cdots x_{9z+6}01(011)^z01(011)^{z+2})^\infty \nonumber\\
											&\preceq(01(011)^{z-2}01(011)^{z}01(011)^{z+2})^\infty \nonumber\\
											&\prec(01(011)^{z-1}01(011)^{z+1}01(011)^z)^\infty\nonumber
										\end{align}
										Hence we have $x_{6z+13}=1$ then $x_{6z+14}=0$ and $x_{6z+15}=1$.
									}
									
								}
								Continue the same argument until $x_{9z+4}=0$, we will get a contraction which implies that $x_{9z+4}=1$ then $x_{9z+5}=0,x_{9z+6}=1$ and $$\sigma^{a+3z+2}((\omega_i))=(01(011)^{2z}0101(011)^z)^\infty.$$
								However we find that $$\sigma^{a+9z+4}((\omega_i))=(0101(011)^{z}01(011)^{2z})^\infty\prec(01(011)^{z-1}01(011)^{z+1}01(011)^z)^\infty$$ which is a contraction to (21).
							\end{itemize}
							Hence we have $x_{3z+5}=1$ then $x_{3z+6}=0,x_{3z+7}=1$ and $$\sigma^a((\omega_i))=(01(011)^{z+1}01x_{3z+8}\cdots x_{9z+6})^\infty.$$
						}
						\item [\ding{194}]
						{
							Continue the same argument as $x_{3z+2}=0$ and $x_{3z+5}=0$ for $x_{3z+8}=0$, we will get a contraction which implies that $x_{3z+8}=1$ and then $x_{3z+9}=0$ and $x_{3z+10}=1$.
						}
						Continue the argument until $x_{9z+2}=0$, we will get a contraction which implies that $x_{9z+2}=1$ then $x_{9z+3}=0,x_{9z+4}=1$ and $$\sigma^a((\omega_i))=(01(011)^{3z}01x_{9z+5}x_{9z+6})^\infty.$$
						\item [\ding{195}]
						{
							If $x_{9z+5}=0$, then we have $x_{9z+6}=1$. \\
							But $\sigma^{a+9z+2}((\omega_i))=(010101(011)^{3z})^\infty\prec(01(011)^{z-1}01(011)^{z+1}01(011)^z)^\infty$ which is a contraction. Hence we have $x_{9z+2}=1$ and then $x_{9z+3}=0$.\\
							However we have $\sigma^{a+9z+5}((\omega_i))=(001(011)^{3z+1})^\infty\prec(01(011)^{z-1}01(011)^{z+1}01(011)^z)^\infty$ which is also a contraction.
						}
					\end{itemize}
					Hence we have $x_3=1$ then $x_4=0,x_5=1$ and $$\sigma^a((\omega_i))=(01101x_6\cdots x_{9z+6})^\infty.$$
					And $\sigma^{a+3}((\omega_i))=(01x_6\cdots x_{9z+6}011)^\infty=:{(y_i)}=(01y_3\cdots y_{9z+3}011)^\infty.$
				}
				\item [2)]
				{
					If $x_6=y_3=0.$\\
					Since $\sigma^{a+3}((\omega_i))=(01y_3\cdots y_{9z+3}011)^\infty\succ(01(011)^{z-1}01(011)^{z+1}01(011)^z)^\infty,$ we have 
					$$(y_i)=(01(011)^{z-1}01y_{3z+2}\cdots y_{9z+3}011)^\infty$$ which is the same situation as the step \ding{192} where $\sigma^a((\omega_i))=(01(011)^{z-1}01x_{3z+2}\cdots x_{9z+6})^\infty$.\\
					Therefore we will get contraction which lead to $$x_6=y_3=1,x_7=y_4=0,x_8=y_5=1.$$
					
				}
				\item [3)]
				{
					Hence we conduct the same argument as \ding{192} and \ding{193} until $x_{9z+3}=y_{9z}=0$ which leads to a contraction which implies that $y_{9z}=1$ then $y_{9z+1}=0,y_{9z+2}=1$ and $$(y_i)=((011)^{3z}01y_{9z+3}011)^\infty.$$
					If $y_{9z+3}=1$, then $(y_i)=((011)^{3z+2})^\infty$ which is of the smallest period 3 which is a contraction. Hence we have $y_{9z+3}=0$.However we find that $$\sigma^{a+9z}(((y_i))=\sigma^{a+9z+3}((\omega_i))=(010011(011)^{3z})^\infty\prec(01(011)^{z-1}01(011)^{z+1}01(011)^z)^\infty.$$ which is also a contraction to (21).
				}
			\end{itemize}
			Hence there is no such (9z+6)-period sequence $(\omega_i)\in \widetilde{K_\beta(t)}$ which is a contraction. Hence we have $$S_\beta(9z+6)\leq [(01(011)^{z-1}01(011)^{z+1}01(011)^z)^\infty]_\beta.$$
			Combined (1) and (2) we get that $$S_\beta(9z+3)= [(01(011)^{z-1}01(011)^{z+1}01(011)^z)^\infty]_\beta,\forall z \in \mathbb{N}_{\geq1}.$$
			
		}
		To sum up, we finally obtain that for $\beta=\phi_3,\ \delta(\beta)=(1^20)^\infty$, we have $S_\beta(1)=0,S_\beta(3)=[(001)^\infty]_\beta,S_\beta(4)=[(0011)^\infty]_\beta,S_\beta(6)=[(001101)^\infty]_\beta$ and
		\begin{equation*}
			S_\beta(p)= \left\{ 
			\begin{aligned}
				&[(01(011)^{m})^\infty]_\beta,p=3m+2,\forall m \in \mathbb{N},\\	&[(01(011)^{z-1}01(011)^z)^\infty]_\beta ,p=3m+1,m=2z,z\in \mathbb{N}_{\geq1},\\
				&[(01(011)^{z-1}01(011)^{z+1})^\infty]_\beta ,p=3m+1,m=2z+1,z\in \mathbb{N}_{\geq1}.\\
				&[(01(011)^{z-1}01(011)^{z-1}01(011)^z)]_\beta^\infty,p=3m,m=3z,z \in \mathbb{N}_{\geq1},\\	&[(01(011)^{z-1}01(011)^z01(011)^z)^\infty]_\beta ,p=3m,m=3z+1,z\in \mathbb{N}_{\geq1},\\
				&[(01(011)^{z-1}01(011)^{z+1}01(011)^z)^\infty]_\beta,p=3m,m=3z+2,z\in \mathbb{N}_{\geq1}.\\
			\end{aligned}  
			\right.        
		\end{equation*}
	\end{proof}
	
\end{them}

\begin{cor}
	By calculation, we obtain that \begin{equation*}
		S_\beta(3m)= \left\{ 
		\begin{aligned}
			&\frac{1+\beta+\beta^{3z+2}+\beta^{6z+1}-\beta^{3z+1}-\beta^{6z}-\beta^{9z-1}-\beta^{9z+1}}{(1-\beta^3)(\beta^{9z}-1)}, m=3z,z \in \mathbb{N}_{\geq1},\\
			&\frac{1+\beta+\beta^{3z+2}+\beta^{6z+4}-\beta^{3z+1}-\beta^{6z+3}-\beta^{9z+2}-\beta^{9z+4}}{(1-\beta^3)(\beta^{9z+3}-1)} ,m=3z+1,z\in \mathbb{N}_{\geq1},\\
			&\frac{1+\beta+\beta^{3z+2}+\beta^{6z+7}-\beta^{3z+1}-\beta^{6z+6}-\beta^{9z+5}-\beta^{9z+7}}{(1-\beta^3)(\beta^{9z+6}-1)},m=3z+2,z\in \mathbb{N}_{\geq1}.\\
		\end{aligned}  
		\right.        
	\end{equation*}
	By take the limit of the expression of $S_\beta(p)$, we get that the value of $S_\beta(p)$ is increasing to $\frac{\beta^{2}+1}{\beta^4-\beta}\approx 0.45626$ as p increasing to infinity.
	
	\begin{figure}[h]  
		\centering  
		\begin{minipage}{.45\textwidth}  
			\centering  
			\includegraphics[width=.8\linewidth]{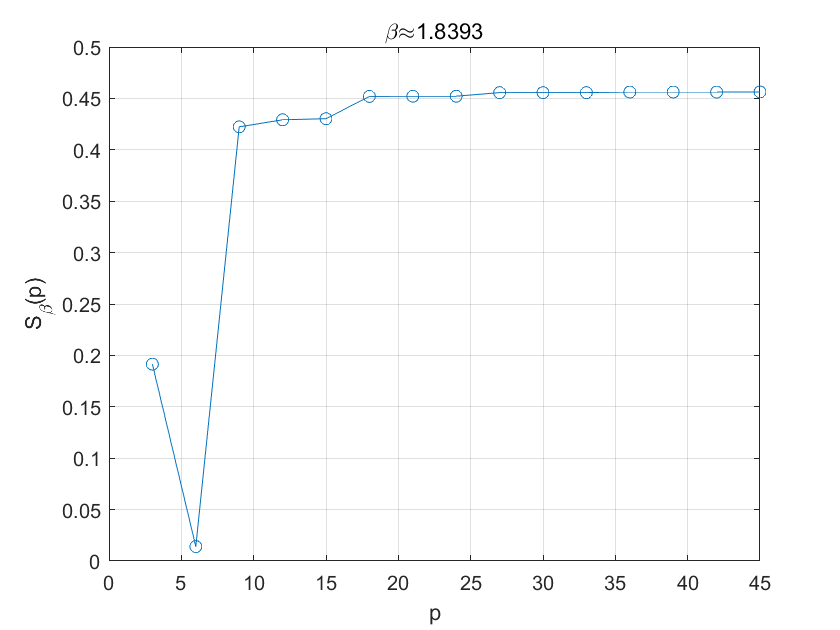}  
			\label{fig:image2}  
		\end{minipage}
		\begin{minipage}{.45\textwidth}  
			\centering  
			\includegraphics[width=.8\linewidth]{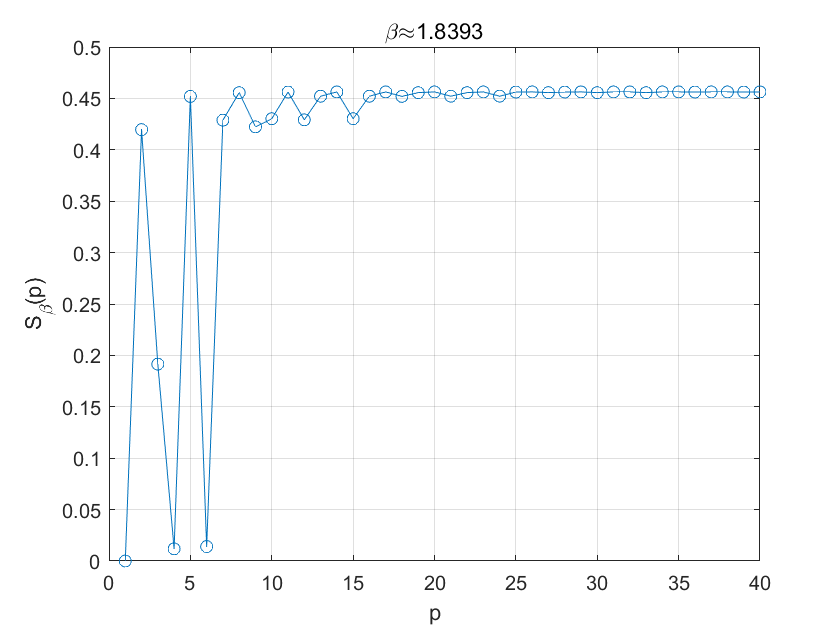}  
			\label{fig:image2}  
		\end{minipage}
		\tiny\caption{Left: the value of $S_\beta(3m)$ for $\beta\approx1.8393$ with the period $\forall m\in\mathbb{N}$ from 1 to 45. Right: the value of $S_\beta(3m)$ for $\beta\approx1.8393$ with the period p from 1 to 40.}
	\end{figure} 		
\end{cor}

\newpage
\addcontentsline{toc}{section}{References}

\end{document}